\documentclass[submission]{eptcs}
\providecommand{\event}{QPL 2013}

\usepackage{latexsym}
\usepackage{amssymb,amsmath,amsthm}
\usepackage{stmaryrd}
\usepackage{xspace}
\usepackage{enumerate}
\usepackage{url}
\usepackage{xypic}
\usepackage{graphicx}
\usepackage{color}
\usepackage{breakurl}

\xyoption{v2}
\xyoption{curve}

\newdimen\proofrulebreadth \proofrulebreadth=.05em
\newdimen\proofdotseparation \proofdotseparation=1.25ex
\newdimen\proofrulebaseline \proofrulebaseline=2ex
\newcount\proofdotnumber \proofdotnumber=3
\let\then\relax
\def\hfi{\hskip0pt plus.0001fil}
\mathchardef\squigto="3A3B
\newif\ifinsideprooftree\insideprooftreefalse
\newif\ifonleftofproofrule\onleftofproofrulefalse
\newif\ifproofdots\proofdotsfalse
\newif\ifdoubleproof\doubleprooffalse
\let\wereinproofbit\relax
\newdimen\shortenproofleft
\newdimen\shortenproofright
\newdimen\proofbelowshift
\newbox\proofabove
\newbox\proofbelow
\newbox\proofrulename
\def\shiftproofbelow{\let\next\relax\afterassignment\setshiftproofbelow\dimen0 }
\def\shiftproofbelowneg{\def\next{\multiply\dimen0 by-1 }%
\afterassignment\setshiftproofbelow\dimen0 }
\def\setshiftproofbelow{\next\proofbelowshift=\dimen0 }
\def\setproofrulebreadth{\proofrulebreadth}

\def\prooftree{
\ifnum  \lastpenalty=1
\then   \unpenalty
\else   \onleftofproofrulefalse
\fi
\ifonleftofproofrule
\else   \ifinsideprooftree
        \then   \hskip.5em plus1fil
        \fi
\fi
\bgroup
\setbox\proofbelow=\hbox{}\setbox\proofrulename=\hbox{}%
\let\justifies\proofover\let\leadsto\proofoverdots\let\Justifies\proofoverdbl
\let\using\proofusing\let\[\prooftree
\ifinsideprooftree\let\]\endprooftree\fi
\proofdotsfalse\doubleprooffalse
\let\thickness\setproofrulebreadth
\let\shiftright\shiftproofbelow \let\shift\shiftproofbelow
\let\shiftleft\shiftproofbelowneg
\let\ifwasinsideprooftree\ifinsideprooftree
\insideprooftreetrue
\setbox\proofabove=\hbox\bgroup$\displaystyle 
\let\wereinproofbit\prooftree
\shortenproofleft=0pt \shortenproofright=0pt \proofbelowshift=0pt
\onleftofproofruletrue\penalty1
}

\def\eproofbit{
\ifx    \wereinproofbit\prooftree
\then   \ifcase \lastpenalty
        \then   \shortenproofright=0pt  
        \or     \unpenalty\hfil         
        \or     \unpenalty\unskip       
        \else   \shortenproofright=0pt  
        \fi
\fi
\global\dimen0=\shortenproofleft
\global\dimen1=\shortenproofright
\global\dimen2=\proofrulebreadth
\global\dimen3=\proofbelowshift
\global\dimen4=\proofdotseparation
\global\count255=\proofdotnumber
$\egroup  
\shortenproofleft=\dimen0
\shortenproofright=\dimen1
\proofrulebreadth=\dimen2
\proofbelowshift=\dimen3
\proofdotseparation=\dimen4
\proofdotnumber=\count255
}

\def\proofover{
\eproofbit 
\setbox\proofbelow=\hbox\bgroup 
\let\wereinproofbit\proofover
$\displaystyle
}%
\def\proofoverdbl{
\eproofbit 
\doubleprooftrue
\setbox\proofbelow=\hbox\bgroup 
\let\wereinproofbit\proofoverdbl
$\displaystyle
}%
\def\proofoverdots{
\eproofbit 
\proofdotstrue
\setbox\proofbelow=\hbox\bgroup 
\let\wereinproofbit\proofoverdots
$\displaystyle
}%
\def\proofusing{
\eproofbit 
\setbox\proofrulename=\hbox\bgroup 
\let\wereinproofbit\proofusing
\kern0.3em$
}

\def\endprooftree{
\eproofbit 
  \dimen5 =0pt
\dimen0=\wd\proofabove \advance\dimen0-\shortenproofleft
\advance\dimen0-\shortenproofright
\dimen1=.5\dimen0 \advance\dimen1-.5\wd\proofbelow
\dimen4=\dimen1
\advance\dimen1\proofbelowshift \advance\dimen4-\proofbelowshift
\ifdim  \dimen1<0pt
\then   \advance\shortenproofleft\dimen1
        \advance\dimen0-\dimen1
        \dimen1=0pt
        \ifdim  \shortenproofleft<0pt
        \then   \setbox\proofabove=\hbox{%
                        \kern-\shortenproofleft\unhbox\proofabove}%
                \shortenproofleft=0pt
        \fi
\fi
\ifdim  \dimen4<0pt
\then   \advance\shortenproofright\dimen4
        \advance\dimen0-\dimen4
        \dimen4=0pt
\fi
\ifdim  \shortenproofright<\wd\proofrulename
\then   \shortenproofright=\wd\proofrulename
\fi
\dimen2=\shortenproofleft \advance\dimen2 by\dimen1
\dimen3=\shortenproofright\advance\dimen3 by\dimen4
\ifproofdots
\then
        \dimen6=\shortenproofleft \advance\dimen6 .5\dimen0
        \setbox1=\vbox to\proofdotseparation{\vss\hbox{$\cdot$}\vss}%
        \setbox0=\hbox{%
                \advance\dimen6-.5\wd1
                \kern\dimen6
                $\vcenter to\proofdotnumber\proofdotseparation
                        {\leaders\box1\vfill}$%
                \unhbox\proofrulename}%
\else   \dimen6=\fontdimen22\the\textfont2 
        \dimen7=\dimen6
        \advance\dimen6by.5\proofrulebreadth
        \advance\dimen7by-.5\proofrulebreadth
        \setbox0=\hbox{%
                \kern\shortenproofleft
                \ifdoubleproof
                \then   \hbox to\dimen0{%
                        $\mathsurround0pt\mathord=\mkern-6mu%
                        \cleaders\hbox{$\mkern-2mu=\mkern-2mu$}\hfill
                        \mkern-6mu\mathord=$}%
                \else   \vrule height\dimen6 depth-\dimen7 width\dimen0
                \fi
                \unhbox\proofrulename}%
        \ht0=\dimen6 \dp0=-\dimen7
\fi
\let\doll\relax
\ifwasinsideprooftree
\then   \let\VBOX\vbox
\else   \ifmmode\else$\let\doll=$\fi
        \let\VBOX\vcenter
\fi
\VBOX   {\baselineskip\proofrulebaseline \lineskip.2ex
        \expandafter\lineskiplimit\ifproofdots0ex\else-0.6ex\fi
        \hbox   spread\dimen5   {\hfi\unhbox\proofabove\hfi}%
        \hbox{\box0}%
        \hbox   {\kern\dimen2 \box\proofbelow}}\doll%
\global\dimen2=\dimen2
\global\dimen3=\dimen3
\egroup 
\ifonleftofproofrule
\then   \shortenproofleft=\dimen2
\fi
\shortenproofright=\dimen3
\onleftofproofrulefalse
\ifinsideprooftree
\then   \hskip.5em plus 1fil \penalty2
\fi
}

Date: Tue, 19 May 1998 16:45:32 +0100
From: Simon Gay <simon@dcs.rhbnc.ac.uk>

I've got another problem when combining
your packages with elsart.cls. The code

doesn't leave enough space below the line in the proof tree, so that
the (a,c) label on the lower arrow runs into the line. It's fine with
article.cls.

\newif\ifignore 

\ignorefalse

\newcommand{\auxproof}[1]{
\ifignore\mbox{}\newline
\textbf{PROOF:} \dotfill\newline
{\it #1}\mbox{}\newline
\textbf{ENDPROOF}\dotfill
\fi}

\newcommand{\C}{\ensuremath{\mathbb{C}}}

\newcommand{\Hil}{\ensuremath{\mathcal{H}}}
\newcommand{\Prob}{\ensuremath{\mathbb{P}}}
\newcommand{\E}{\ensuremath{\mathbb{E}}}

\newcommand{\Proj}{\ensuremath{\mathrm{Proj}}}

\newcommand{\D}{\ensuremath{\mathcal{D}}}
\newcommand{\G}{\ensuremath{\mathcal{G}}}

\newcommand{\F}{\ensuremath{\mathcal{F}}}
\newcommand{\Ctr}{\ensuremath{\mathcal{Z}}}
\newcommand{\Rdn}{\ensuremath{\mathcal{R}}}

\newcommand{\klcomp}{\ensuremath{\mathrel{\raisebox{.15em}{$\scriptscriptstyle\odot$}}}}

\newcommand{\ch}{\ensuremath{\mathrm{char}}}

\newcommand{\intd}{{\kern.2em}\mathrm{d}{\kern.03em}}

\newcommand{\isqrt}[1]{\frac{1}{\sqrt{#1}}}

\newcommand{\bra}[1]{\ensuremath{\langle #1 |}}
\newcommand{\ket}[1]{\ensuremath{| #1 \rangle}}
\newcommand{\braket}[2]{\ensuremath{\langle #1 | #2 \rangle}}
\newcommand{\braketI}[1]{\braket{#1}{#1}}
\newcommand{\ketbra}[2]{\ensuremath{ | #1 \rangle \langle #2 | }}
\newcommand{\ketbraI}[1]{\ketbra{#1}{#1}}
\newcommand{\given}[2]{\ensuremath{#1{\kern.1em}|{\kern.1em}#2}}
\newcommand{\joint}[2]{\ensuremath{#1 \wedge #2}}

\newcommand{\after}{\mathrel{\circ}}

\newcommand{\setin}[3]{\{#1\in#2\;|\;#3\}}

\newcommand{\allin}[3]{\forall_{#1\in#2}.\,#3}

\newcommand{\tuple}[1]{\langle#1\rangle}
\newcommand{\NNO}{\mathbb{N}}

\newcommand{\op}[1]{#1^{\textrm{op}}}

\newcommand{\charac}{\ensuremath{\textsl{char}}}
\newcommand{\graph}{\ensuremath{\textsl{gr}}}

\newcommand{\tr}{\ensuremath{\textsl{tr}}}

\newcommand{\sotimes}{\mathrel{\raisebox{.05pc}{$\scriptstyle \otimes$}}}

\newcommand{\Sets}{\ensuremath{\textbf{Sets}}\xspace}

\newcommand{\EA}{\ensuremath{\textbf{EA}}\xspace}

\newcommand{\CHaus}{\ensuremath{\textbf{CHaus}}\xspace}
\newcommand{\EMod}{\ensuremath{\textbf{EMod}}\xspace}

\newcommand{\Meas}{\ensuremath{\textbf{Meas}}\xspace}

\newcommand{\CstarMap}[1]{\ensuremath{\textbf{Cstar}_{#1}}\xspace}
\newcommand{\CstarPU}{\CstarMap{\textrm{PU}}}
\newcommand{\CstarCPU}{\CstarMap{\textrm{cPU}}}
\newcommand{\CstarMIU}{\CstarMap{\textrm{MIU}}}

\newcommand{\CCstarMap}[1]{\ensuremath{\textbf{CCstar}_{#1}}\xspace}
\newcommand{\CCstarPU}{\CCstarMap{\textrm{PU}}}

\newcommand{\Pred}{\ensuremath{\textsl{Pred}}\xspace}

\newcommand{\scalar}{\mathrel{\bullet}}
\newcommand{\inv}{\mathop{\rlap{\raisebox{.3ex}{${\kern.6ex}\cdot$}}-}}
\newcommand{\id}{\ensuremath{\mathrm{id}}}
\newcommand{\idmap}[1][]{\ensuremath{\mathrm{id}_{#1}}}

\newcommand{\subst}[1]{\ensuremath{#1^{\sharp}}}

\newcommand{\Kl}{\mathcal{K}{\kern-.2ex}\ell}
\newcommand{\KlN}{\mathcal{K}{\kern-.2ex}\ell_{\NNO}}
\newcommand{\EM}{\mathcal{E}{\kern-.2ex}\mathcal{M}}
\newcommand{\Ef}{\ensuremath{\mathcal{E}{\kern-.5ex}f}}
\newcommand{\DM}{\ensuremath{\mathcal{D}{\kern-.85ex}\mathcal{M}}}
\newcommand{\FP}{\mathcal{F}{\kern-.2ex}\mathcal{P}}
\newcommand{\CR}{\ensuremath{\mathcal{C}_{\mathcal{R}}}}

\newcommand{\QEDbox}{\square}
\newcommand{\QED}{\hspace*{\fill}$\QEDbox$}

\newtheorem{theorem}{Theorem}

\newtheorem{lemma}[theorem]{Lemma}
\newtheorem{corollary}[theorem]{Corollary}
\newtheorem{example}[theorem]{Example}

\newenvironment{myproof}[1][Proof]%
   { \begin{trivlist}%
     \item[\hskip \labelsep {\bfseries #1}]%
   }%
   { \end{trivlist}%
   }

\def\authorrunning{Furber and Jacobs}
\def\titlerunning{Towards a Categorical Account of Conditional Probability}

\title{Towards a Categorical Account of \\ Conditional Probability}

\author{Robert Furber and Bart Jacobs
\institute{
Institute for Computing and Information Sciences (iCIS), \\
Radboud University Nijmegen, The Netherlands. \\
Web addresses: \url{www.cs.ru.nl/~rfurber} and 
   \url{www.cs.ru.nl/~bart}\\[.5em]
}
}

\begin{document}
\maketitle

\begin{abstract}
This paper presents a categorical account of conditional probability,
covering both the classical and the quantum case. Classical
conditional probabilities are expressed as a certain
``triangle-fill-in'' condition, connecting marginal and joint
probabilities, in the Kleisli category of the distribution monad. The
conditional probabilities are induced by a map together with a
predicate (the condition). The latter is a predicate in the logic of
effect modules on this Kleisli category.

This same approach can be transferred to the category of
$C^*$-algebras (with positive unital maps), whose predicate logic is
also expressed in terms of effect modules. Conditional probabilities
can again be expressed via a triangle-fill-in property. In the literature, 
there are several proposals for what quantum conditional probability should
be, and also there are extra difficulties not present in the classical case.
At this stage, we only describe quantum systems with classical parametrization.

\end{abstract}

\renewcommand{\arraystretch}{1.3}
\setlength{\arraycolsep}{2pt}

\section{Introduction}

In the categorical description of probability theory, several monads
play an important role. The main ones are the discrete probability
monad $\D$ on the category $\Sets$ of sets and functions, and the Giry
monad $\G$, for continuous probability, on the category $\Meas$ of
measurable spaces and measurable functions. The Kleisli categories of
these monads have suitable probabilistic matrices as morphisms, which
capture probabilistic transition systems (and Markov
chains). Additionally, more recent monads of interest are the
expectation monad~\cite{JacobsM12b} and the Radon
monad~\cite{FurberJ13a}.

The first contribution of this paper is a categorical reformulation
of classical (discrete) conditional probability as a ``triangle-fill in''
property in the Kleisli category $\Kl(\D)$ of the distribution monad.
Abstractly, this fill-in property appears as follows.
\begin{equation}
\label{CondProbTrianglePattern}
\vcenter{\xymatrix@R+1pc@C+1pc{
& X \ar[dl]_{\mbox{marginal probability}\quad}
   \ar[dr]^{\qquad\mbox{joint probability}} & \\
X + X \ar@{-->}[rr]_-{\begin{array}{c} \mbox{conditional} \\[-.6em]
   \mbox{probability} \end{array}} & & Y + Y
}}
\end{equation}

\noindent This diagram incorporates the idea that `conditional'
$\times$ `marginal' = `joint'. This idea is illustrated in two
examples: first in the simpler non-parametrized case, and later also
in parametrized form.

The same idea can be expressed in other Kleisli categories, of the
other monads mentioned above. But a more challenging issue is to
transfer this approach to the quantum case. This constitutes the main
part (and contribution) of this paper. We interpret the above triangle
in the opposite of the category of $C^*$-algebras, with positive
unital maps, using effects as predicates. In this quantum case the
situation becomes more subtle, and at this preliminary stage of
investigation we only present a non-parametrized example, namely the
Elitzur-Vaidman bomb tester~\cite{Elitzur93}.

Our work relates to the pre-existing literature as follows. Bub \cite{bub07} interprets the projection postulate during a measurement as an instance of Bayesian updating of the quantum state. His formulas \cite[(21) and (22)]{bub07}  for the special case that $B = B(\Hil)$ (and the state is normal, which is satisfied automatically if $\dim \Hil < \infty$) agree with our formula \eqref{UnParamGivenEqn}. 

We can see this as follows. Bub has, for $a,b \in \Proj(\Hil)$:
$$\begin{array}{rclcrcl}
\Prob_\rho(b | a)
& = &
\tr(\rho' b)
& \qquad\mbox{where}\qquad &
\rho'
& = &
\frac{a\rho a}{\tr(a \rho a)}.
\end{array}$$

\noindent This can be rearranged:
$$\begin{array}{rcccccccl}
\Prob_\rho(b | a) 
& = & 
\displaystyle\tr \left(\frac{a \rho a}{\tr(a \rho a)} b \right)
& = & 
\displaystyle\frac{\tr(a \rho a b)}{\tr(a \rho a)}
& = & 
\displaystyle\frac{\tr(\rho (a b a))}{\tr(\rho a^2)}
& = & 
\displaystyle\frac{\tr(\rho (a b a))}{\tr(\rho a)}.
\end{array}$$

\noindent If we reinterpret the $\rho$s as maps $B(\Hil) \rightarrow
\C$, i.e. normal states, we get:
\[
\Prob_\rho(b | a) = \frac{\rho(a b a)}{\rho(a)}
\]
\noindent We now see this agrees with \eqref{UnParamGivenEqn}.

There is also the quantum conditional probability definition of Leifer
and Spekkens~\cite{leifer12} --- expressed graphically
in~\cite{coecke12}. This work is based on the probabilistic case
where, instead of being expressed in terms of probabilities of
predicates, the conditional probability is formulated
in~\cite{leifer12} using a random variable that completely determines
the elements of the underlying probability space. This seems to lead
to a different formula from ours, but precise comparison is left to
future work.




\section{Discrete probability, categorically}\label{DiscProbSec}

To describe finite discrete probabilities categorically one uses the
distribution monad $\D\colon\Sets \rightarrow \Sets$. It maps a set
$X$ to the set $\D(X)$ of probability distributions over $X$, which we
describe as formal finite convex sums:
$$\textstyle \sum_{i} r_{i}\ket{x_i}
\qquad\mbox{where}\qquad
x_{i}\in X \mbox{ and } r_{i}\in [0,1] \mbox{ satisfy } \sum_{i}r_{i}=1.$$

\noindent We use the ``ket'' notation $\ket{-}$ to distinguish
elements $x\in X$ and their occurrences in formal sums. Each function
$f\colon X \rightarrow Y$ gives a function $\D(f) \colon \D(X)
\rightarrow \D(Y)$, where:
$$\begin{array}{rcl}
\D(f)\big(\sum_{i}r_{i}\ket{x_i}\big)
& = &
\sum_{i}r_{i}\ket{f(x_{i})}.
\end{array}$$

\noindent The unit $\eta \colon X \rightarrow \D(X)$ of this
distribution monad $\D$ sends $x\in X$ to the singleton/Dirac
distribution $\eta(x) = 1\ket{x}$. The multiplication $\mu \colon
\D^{2}(X) \rightarrow \D(X)$ is given by:
$$\begin{array}{rclcrcl}
\mu\big(\sum_{i} r_{i}\ket{\varphi_i}\big)
& = &
\sum_{i,j} (r_{i}s_{ij})\ket{x_{ij}}
& \qquad\mbox{if}\qquad &
\varphi_{i}
& = &
\sum_{j} s_{ij}\ket{x_{ij}}.
\end{array}$$

Like for any monad, one can form the Kleisli category $\Kl(\D)$. In
this case we get the category of sets and stochastic matrices, as
the objects of $\Kl(\D)$ are sets, and its maps
$X\rightarrow Y$ are functions $X \rightarrow \D(Y)$. The unit
function $\eta\colon X \rightarrow \D(X)$ is then the identity map $X
\rightarrow X$ in $\Kl(\D)$. Composition of $f\colon X \rightarrow Y$
and $g\colon Y \rightarrow Z$ in $\Kl(\D)$ yields a map $g \klcomp f
\colon X \rightarrow Z$, which, as a function $X \rightarrow
\D(Z)$ is given by $g \klcomp f = \mu \after \D(g) \after
f$. Explicitly:
$$\begin{array}{rclcrclcrcl}
(g \klcomp f)(x)
& = &
\sum_{i,j} (r_{i}s_{ij})\ket{z_{ij}}
& \quad\mbox{if}\quad &
f(x) 
& = &
\sum_{i}r_{i}\ket{y_i}
& \quad\mbox{and}\quad &
g(y_{i})
& = &
\sum_{j}s_{ij}\ket{z_{ij}}.
\end{array}$$

\noindent There is a forgetful functor $\Kl(\D) \rightarrow \Sets$,
sending $X$ to $\D(X)$ and $f$ to $\mu \after \D(f)$. It has a left
adjoint $\mathcal{F} \colon \Sets \rightarrow \Kl(\D)$ which is the
identity on objects and sends $f$ to $\eta\after f$.

Products and coproducts of sets, with their projections $\pi_i$ and
coprojections $\kappa_i$ are written as:
$$\xymatrix{
X & X\times Y\ar[l]_-{\pi_1}\ar[r]^-{\pi_2} & Y
& &
X\ar[r]^-{\kappa_{1}} & X+Y & Y\ar[l]_-{\kappa_2}
}$$

\noindent There are associated tuples $\tuple{f,g} \colon Z
\rightarrow X\times Y$ and cotuples $[h,k]\colon X+Y\rightarrow
Z$. The empty product is a singleton set, typically written as $1$,
and the empty coproduct is the empty set $0$.

The category $\Kl(\D)$ inherits these coproducts $(+,0)$ from $\Sets$,
with coprojections $\mathcal{F}(\kappa_{i}) = \eta \after \kappa_{i}$,
and cotupling $[f,g]$ as in $\Sets$. The products $(\times, 1)$ from
$\Sets$ form a tensor product --- not a cartesian product --- on
$\Kl(\D)$; hence we write $\otimes$ in $\Kl(\D)$ for $\times$. But
because the tensor unit $1$ is also final in $\Kl(\D)$, since $\D(1)
\cong 1$, we have a tensor with projections in $\Kl(\D)$. We shall
write $\pi_{i} \colon X_{1}\otimes X_{2} \rightarrow X_{i}$ for the
resulting projections in $\Kl(\D)$, which are functions
$\mathcal{F}(\pi_{i}) = \eta \after \pi_{i} \colon X_{1}\times X_{2}
\rightarrow \D(X_{i})$. This forms the background for the following
result. It uses \emph{marginals}, which, for a Kleisli map $f\colon X
\rightarrow Y_{1}\otimes Y_{2}$ are obtained by post-composition
$\pi_{i} \klcomp f = \D(\pi_{i}) \after f \colon X
\rightarrow Y_{i}$.

\auxproof{
Proof that $\pi_i \klcomp f = \D(\pi_i) \circ f$:
\begin{eqnarray*}
\pi_i \klcomp f & = & \mathcal{F}(\pi_i) \klcomp f \\
 & = & (\eta_{Y_i} \circ \pi_i) \klcomp f \\
 & = & \mu_{Y_i} \circ \D(\eta_{Y_i} \circ \pi_i) \circ f \\
 & = & \mu_{Y_i} \circ \D(\eta_{Y_i}) \circ \D(\pi_i) \circ f \\
 & = & \D(\pi_i) \circ f
\end{eqnarray*}
}

\begin{lemma}
\label{KlDGraphLem}
In $\Kl(\D)$ there is a bijective correspondence:
$$\begin{prooftree}
{\xymatrix{ X \ar[r]^-{f} & Y }}
\Justifies
{\xymatrix{ X \ar[r]_-{g} & X\otimes Y \mbox{ with } 
   \pi_{1} \klcomp g = \idmap[X]}}
\end{prooftree}$$
\end{lemma}

\begin{myproof}
The condition $\pi_{1} \klcomp g = \idmap$ means that if
$g(x) = \sum_{i}r_{i}\ket{(x_{i},y_{i})}$, then $x_{i} = x$ for all
$i$. Hence $g$ corresponds to a function $X \rightarrow \D(Y)$. \QED
\end{myproof}

Below we shall write $\graph(f)\colon X \rightarrow \D(X\times Y)$ for
this ``graph'' map corresponding to $f\colon X \rightarrow \D(Y)$,
where, explicitly,
\begin{equation}
\label{KlDGraphEqn}
\begin{array}{rclcrcl}
\graph(f)(x)
& = &
\sum_{i}r_{i}\ket{(x,y_{i})}
& \quad\mbox{if}\quad &
f(x)
& = &
\sum_{i}r_{i}\ket{y_i}.
\end{array}
\end{equation}

Now that we have a category $\Kl(\D)$ to model probabilistic
transitions, we add a logic to it. Categorically this takes the form of
a functor, or indexed category, $\Pred \colon \Kl(\D) \rightarrow
\op{\EMod}$, where $\EMod$ is the category of effect modules (see
\textit{e.g.}~\cite{JacobsM12c}). We briefly explain the relevant
definitions. 

To start, let $M = (M, \ovee, 0)$ be a partial commutative monoid,
where $\ovee$ is a partial operation $M\times M \rightarrow M$ that is
commutative and associative, in a suitable sense, and has $0$ has unit
element. One can think of the unit interval $[0,1]$ with addition $+$
and $0$. Such an $M$ is called an \emph{effect algebra} if there is a
unary ``orthocomplement'' operation $(-)^{\perp} \colon M\rightarrow
M$ satisfying both:
\begin{itemize}
\item $x^{\perp}\in E$ is the unique element in $E$ with $x\ovee
  x^{\perp} = 1$, where $1 = 0^\perp$;

\item $x\ovee 1$ is defined only when $x=0$.
\end{itemize}

\noindent On $[0,1]$ on has $r^{\perp} = 1 - r$ as orthocomplement.

A morphism of effect algebras $f\colon M \rightarrow N$ is a function
between the underlying sets satisfying $f(1) = 1$ and: if $x\ovee y$
is defined, then so is $f(x)\ovee f(y)$, and $f(x\ovee y) = f(x) \ovee
f(y)$. This yields a category which we write as $\EA$.

An \emph{effect module} is an effect algebra $M$ with a (total) scalar
multiplication $r\scalar x \in M$ for $r\in[0,1]$ and $x\in M$,
preserving $\ovee$ in both coordinates separately, satisfying
$1\scalar x = x$, and $r\scalar (s\scalar x) = (r \cdot s) \scalar x$.
A map of effect modules is a map of effect algebras that preserves 
the scalar multiplication. This yields a category $\EMod$.

For a set $X\in\Kl(\D)$ we define $\Pred(X) = [0,1]^{X}$, the set of
fuzzy predicates on $X$. There are a true and false predicates, $1 =
x \mapsto 1$ and $0 = x \mapsto 0$. For two fuzzy predicates $p,q\in
    [0,1]^{X}$ a sum $p \ovee q\in [0,1]^{X}$ exists if $p(x)+q(x)\leq
    1$, for all $x\in X$; then $(p\ovee q)(x) = p(x) + q(x)$. Further,
    there is an orthocomplement operation $p^{\perp}(x) =
    1-p(x)$. One has, for instance, $p^{\perp\perp} = p$ and
    $p \ovee p^{\perp} = 1$. There is also a scalar multiplication on
    fuzzy predicates: for $r\in[0,1]$ one defines $(r\scalar p)(x) =
    r\cdot p(x)$.

Each Kleisli map $f\colon X \rightarrow Y$ yields a functor $\subst{f} =
\Pred(f) \colon \Pred(Y) \rightarrow \Pred(X)$, which is commonly
called substitution; it is given by:
\begin{equation}
\label{SubstProbEqn}
\begin{array}{rclcrcl}
\subst{f}(q)(x)
& = &
\sum_{i} r_{i}\cdot q(y_{i})
& \qquad\mbox{if}\qquad &
f(x)
& = &
\sum_{i}r_{i}\ket{y_{i}}.
\end{array}
\end{equation}

\noindent This $\subst{f}$ is a map of effect modules $[0,1]^{Y}
\rightarrow [0,1]^{X}$.

For each set $X$ there is a special predicate $\Omega_{X}
\in \Pred(X+X) = [0,1]^{X+X}$, namely $\Omega(\kappa_{1}x) = 1$ and
$\Omega(\kappa_{2}x) = 0$. For each predicate $p \in [0,1]^{X}$ there is a
\emph{characteristic} map $\charac_{p} \colon X \rightarrow X+X$
in $\Kl(\D)$ with $\subst{\charac_{p}}(\Omega) = p$. This characteristic
map is defined as convex sum:
\begin{equation}
\label{CharacProbEqn}
\begin{array}{rcccl}
\charac_{p}(x)
& = &
p(x)\ket{\kappa_{1}x} + p^{\perp}(x)\ket{\kappa_{2}x}
& = &
p(x)\ket{\kappa_{1}x} + (1-p(x))\ket{\kappa_{2}x}.
\end{array}
\end{equation}

\noindent These characteristic maps play an important role below, and are further discussed in \cite{Jacobs13b}.

\section{Conditional discrete probability}\label{CondDiscProbSec}

This section reviews classical conditional probability, in the
discrete case. A simple example is first described in standard
terminology, and then reformulated in categorical form, by using the
fuzzy predicate logic $\Pred\colon \Kl(\D) \rightarrow\nobreak \op{\EMod}$
over the Kleisli category of the (discrete) probability monad $\D$.
The example is extended to ``parametrized'' form, and again formulated
in categorical terms.

\begin{example}
\label{HairEx}
In this first illustration we describe a simple situation, involving a
set of genders $G = \{M, W\}$ with a distribution $f =
\frac{2}{3}\ket{M} + \frac{1}{3}\ket{W}$ of men and women. Assume that
the probability of having long hair is $\frac{3}{10}$ for men and
$\frac{8}{10}$ for women. More formally this is written as
$\Prob[\given{\ell}{M}] = \frac{3}{10}$ and $\Prob[\given{\ell}{W}] =
\frac{8}{10}$, where $\ell$ stands for `long hair'. We now ask
ourselves the typical conditional probability question: suppose we see
someone with long hair, what is the probability that the person is a
man/woman?

One then proceeds as follows. The joint probabilities are given by:
$$\begin{array}{rccclcrcccl}
\Prob[\joint{M}{\ell}]
& = &
\frac{2}{3}\cdot \frac{3}{10}
& = &
\frac{1}{5}
& \qquad &
\Prob[\joint{W}{\ell}]
& = &
\frac{1}{3}\cdot \frac{8}{10}
& = &
\frac{4}{15}.
\end{array}$$

\noindent And the marginal probability of seeing long hair is:
$$\begin{array}{rcccccl}
\Prob[\ell]
& = &
\Prob[\joint{M}{\ell}] + \Prob[\joint{W}{\ell}]
& = &
\frac{1}{5} + \frac{4}{15}
& = &
\frac{7}{15}.
\end{array}$$ 

\noindent We then obtain the required conditional probabilities:
$$\begin{array}{rccccclcrcccccl}
\Prob[\given{M}{\ell}]
& = &
\displaystyle\frac{\Prob[\joint{M}{\ell}]}{\Prob[\ell]}
& = &
\displaystyle\frac{\;\frac{1}{5}\;}{\frac{7}{15}}
& = &
\frac{3}{7}
& \qquad &
\Prob[\given{W}{\ell}]
& = &
\displaystyle\frac{\Prob[\joint{W}{\ell}]}{\Prob[\ell]}
& = &
\displaystyle\frac{\;\frac{4}{15}\;}{\frac{7}{15}}
& = &
\frac{4}{7}.
\end{array}$$

\noindent By construction we have ``conditional $\,\cdot\,$ marginal =
joint'' since $\Prob[\given{M}{\ell}] \cdot \Prob[\ell] =
\Prob[\joint{M}{\ell}]$, as suggested
in~\eqref{CondProbTrianglePattern}.  It will be elaborated below.
\end{example}

We now reformulate this example in categorical form. The distribution
$f = \frac{2}{3}\ket{M} + \frac{1}{3}\ket{W}$ corresponds to a map
$f\colon 1 \rightarrow G$ in the Kleisli category $\Kl(\D)$, where $1
= \{0\}$ is the final (singleton) set and $G = \{M, W\}$ is the
two-element set of genders. In this correspondence we identify $f$
with the value $f(0)\in \D(G)$, for the sole element $0\in 1$. The
likelihood of having long hair corresponds to a fuzzy predicate
$\ell\in\Pred(G) = [0,1]^{G}$ on the set $G$, given by $\ell(M) =
\frac{3}{10}$, $\ell(W) = \frac{8}{10}$. The associated characteristic
map $\charac_{\ell} \colon G \rightarrow G+G$ in $\Kl(\D)$ is, according
to~\eqref{CharacProbEqn}:
$$\begin{array}{rclcrcl}
\charac_{\ell}(M)
& = &
\frac{3}{10}\ket{\kappa_{1}M} + \frac{7}{10}\ket{\kappa_{2}M}
& \qquad &
\charac_{\ell}(W)
& = &
\frac{8}{10}\ket{\kappa_{1}W} + \frac{2}{10}\ket{\kappa_{2}W}.
\end{array}$$

\noindent The left $\kappa_1$-option in the coproduct $G+G$ thus
captures the probability that the predicate (in this case $\ell$) is
true, and the right $\kappa_{2}$-option is for false. The composite
map $\charac_{\ell} \klcomp f \colon 1 \rightarrow G+G$ in the Kleisli
category now describes the \emph{joint} probability:
$$\begin{array}{rcl}
\charac_{\ell} \klcomp f
& = &
\sum_{z\in G+G} \big(\sum_{g\in G} f(g)\cdot \charac_{\ell}(g)(z)\big) \, \ket{z} \\
& = &
\frac{2}{3}\cdot\frac{3}{10}\ket{\kappa_{1}M} +
   \frac{2}{3}\cdot\frac{7}{10}\ket{\kappa_{2}M} +
   \frac{1}{3}\cdot\frac{8}{10}\ket{\kappa_{1}W} +
   \frac{1}{3}\cdot\frac{2}{10}\ket{\kappa_{2}W} \\
& = &
\frac{1}{5}\ket{\kappa_{1}M} +
   \frac{7}{15}\ket{\kappa_{2}M} +
   \frac{4}{15}\ket{\kappa_{1}W} +
   \frac{1}{15}\ket{\kappa_{2}W} \\
& = &
\Prob[\joint{M}{\ell}]\ket{\kappa_{1}M} +
    \Prob[\joint{M}{\ell^{\perp}}]\ket{\kappa_{2}M} +
    \Prob[\joint{W}{\ell}]\ket{\kappa_{1}W} +
    \Prob[\joint{W}{\ell^{\perp}}]\ket{\kappa_{2}W}.
\end{array}$$

\noindent The substituted predicate $\subst{f}(\ell) \in \Pred(1) =
          [0,1]^{1} \cong [0,1]$, defined in~\eqref{SubstProbEqn},
          gives the \emph{marginal} probability $\Pr[\ell] \in [0,1]$:
$$\begin{array}{rcccccccl}
\subst{f}(\ell)
& = &
\sum_{g\in G}f(g)\cdot\ell(g)
& = &
f(M)\cdot\ell(M) + f(W)\cdot\ell(W) 
& = &
\frac{2}{3}\cdot \frac{3}{10} + \frac{1}{3}\cdot \frac{8}{10} 
& = &
\frac{7}{15}.
\end{array}$$


The conditional probabilities can be organized into two
maps $f|\ell, f|\ell^{\perp} \colon 1 \rightarrow G$ in $\Kl(\D)$, namely:
$$\begin{array}{rcl}
\given{f}{\ell}
& = &
\frac{1}{\Prob[\ell]}\big(\Prob[\joint{M}{\ell}]\,\ket{M} + 
   \Prob[\joint{W}{\ell}]\, \ket{W}\big)
\hspace*{\arraycolsep}=\hspace*{\arraycolsep}
\frac{15}{7}\big(\frac{1}{5}\ket{M} + \frac{4}{15}\ket{W}\big)
\hspace*{\arraycolsep}=\hspace*{\arraycolsep}
\frac{3}{7}\ket{M} + \frac{4}{7}\ket{W} \\
f|\ell^{\perp}
& = &
\frac{1}{\Prob[\ell^{\perp}]}\big(\Prob[\joint{M}{\ell^{\perp}}]\,\ket{M} + 
   \Prob[\joint{W}{\ell^{\perp}}]\, \ket{W}\big) 
\hspace*{\arraycolsep}=\hspace*{\arraycolsep}
\frac{15}{8}\big(\frac{7}{15}\ket{M} + \frac{1}{15}\ket{W}\big) 
\hspace*{\arraycolsep}=\hspace*{\arraycolsep}
\frac{7}{8}\ket{M} + \frac{1}{8}\ket{W}.
\end{array}$$

\noindent The first distribution $\given{f}{\ell}$ gives the
probabilities for men and women under the assumption that you see long
hair; similarly, $\given{f}{\ell^{\perp}}$ gives these probabilities
if you do not see long hair.

The final observation is that these two maps $\given{f}{\ell}$ and
$\given{f}{\ell^\perp}$ make the following triangle in the Kleisli
category $\Kl(\D)$ commute, like in
pattern~\eqref{CondProbTrianglePattern}:
\begin{equation}
\label{HairTriangle}
\vcenter{\xymatrix@R+1pc@C+1pc{
& 1\ar[dl]_{\charac_{\subst{f}(\ell)}}
   \ar[dr]^{\charac_{\ell} \klcomp f} & \\
1+1\ar[rr]_-{(\given{f}{\ell}) + (\given{f}{\ell^{\perp}})} & & G+G
}}
\end{equation}

This simple hair example is ``non-parametrized'', in the sense that
in the above triangle we have the final/singleton set $1$ at the
top. More generally, we can start with a Kleisli map $f\colon X
\rightarrow Y$ and predicate on $X\otimes Y$.

\begin{example}
\label{CountryHairEx}
Suppose we now have two different countries $A,B$ which have different
gender distributions and different distributions of long and short
hair. We will use $C = \{ A, B \}$ as the set of countries, with given
gender distributions captured by a Kleisli map $f\colon C \rightarrow
\D(G)$, where $G = \{M, W\}$ is the set of genders like in
Example~\ref{HairEx}:
$$\begin{array}{rclcrcl}
f(A) & = & \frac{9}{20}\ket{M} + \frac{11}{20}\ket{W}
& \qquad &
f(B) & = & \frac{1}{2}\ket{M} + \frac{1}{2}\ket{W}.
\end{array}$$

\noindent The probabilities of having long hair depend on both $C$ and
$G$ and are already given in some way, formalised via predicate $L \in
\Pred(C \otimes G)$ with:
$$\begin{array}{rclcrclcrclcrcl}
L(A,M) & = & \frac{1}{10}
& \quad &
L(B,M) & = & \frac{2}{10}
& \quad &
L(A,W) & = & \frac{8}{10}
& \quad &
L(B,W) & = & \frac{9}{10}.
\end{array}$$

\noindent When instead of a proper set $C$ we had a trivial
(singleton) set $1$, getting the joint probability distribution was a
simple matter of composition in $\Kl(\D)$. However, we now have
$f\colon C \rightarrow G$ but $\ch_L \colon C \otimes G \rightarrow C
\otimes G + C \otimes G$, so that is out of the question.  To solve this
we will define a map $j = \joint{f}{L} \colon C \rightarrow G + G$ as
the composite of the following maps, using the definition of $\graph$
from equation~\eqref{KlDGraphEqn}:
\begin{equation}
\label{jointProb}
\vcenter{\xymatrix{
C \ar[rr]^-{j = \joint{f}{L}} \ar[d]_{\graph(f)} & & G + G \\
C \otimes G \ar[rr]_-{\ch_L} & &
   C \otimes G + C \otimes G \ar[u]_{\F(\pi_2 + \pi_2)}
}}
\end{equation}

\noindent This produces the correct joint probability. First the case
of $A \in C$:
$$\begin{array}{rcl}
j(A) 
& = & 
\big(\D(\pi_2 + \pi_2) \after (\ch_{L} \klcomp \graph(f))\big)(A) \\
& = & 
\big(\D(\pi_2 + \pi_2) \after \ch_{L}^{\$}\big)
   \big(\sum_{g \in G}f(A)(g)\ket{A,g}\big) \\
& = & 
\D(\pi_2 + \pi_2)\big(\sum_{g\in G} 
   L(A,g)\cdot f(A)(g)\ket{\kappa_{1}(A,g)} +
   L^{\perp}(A,g)\cdot f(A)(g)\ket{\kappa_{2}(A,g)}\big) \\
& = &
\D(\pi_2 + \pi_2)\big(\frac{1}{10}\cdot\frac{9}{20}\ket{\kappa_1(A,M)} + 
   \frac{9}{10}\cdot \frac{9}{20}\ket{\kappa_2(A,M)} \\
& & \hspace*{5em} \;+\; 
   \frac{8}{10}\cdot \frac{11}{20}\ket{\kappa_1(A,W)} + 
   \frac{2}{10}\cdot\frac{11}{20}\ket{\kappa_2(A,W)}\big) \\
& = & 
\frac{9}{200}\ket{\kappa_1M} + \frac{81}{200}\ket{\kappa_2M} + 
   \frac{88}{200}\ket{\kappa_1W} + \frac{22}{200}\ket{\kappa_2W}.
\end{array}$$

\noindent Similarly one obtains the distribution $j(B)\in \D(G+G)$, namely:
$$\begin{array}{rcl}
j(B)
& = &
\frac{2}{20}\ket{\kappa_1M} + \frac{8}{20}\ket{\kappa_2M} + 
   \frac{9}{20}\ket{\kappa_1W} + \frac{1}{20}\ket{\kappa_2W}.
\end{array}$$

\auxproof{
Now the case of $B \in C$:
\begin{eqnarray*}
j(B) & = & (\pi_2 + \pi_2)(\ch_L(\graph(f)(B))) \\
 & = & (\pi_2 + \pi_2)\left(\ch_L\left(\sum_{g \in G}f(B)(g)\ket{B,g}\right)\right) \\
 & = & (\pi_2 + \pi_2)(\ch_L(\frac{1}{2}\ket{B,M} + \frac{1}{2}\ket{B,W})) \\
 & = & (\pi_2 + \pi_2)(\frac{1}{2}\frac{2}{10}\ket{\kappa_1(B,M)} + \frac{1}{2}\frac{8}{10}\ket{\kappa_2(B,M)} + \frac{1}{2}\frac{9}{10}\ket{\kappa_1(B,W)} + \frac{1}{2}\frac{1}{10}\ket{\kappa_2(B,W)} \\
 & = & 
\frac{2}{20}\ket{\kappa_1M} + \frac{8}{20}\ket{\kappa_2M} + 
   \frac{9}{20}\ket{\kappa_1W} + \frac{1}{20}\ket{\kappa_2W}
\end{eqnarray*}
}

We now calculate the marginal probability of $L$, getting rid of the
dependence on $G$. To do this, we compute $\subst{\graph(f)}(L) \in \Pred(C)$ as:
$$\begin{array}{rcl}
\subst{\graph(f)}(L)(A) 
& = & 
\sum_{(c,g) \in C \times G} \graph(f)(A)(c,g) \cdot L(c,g) \\
& = & 
\frac{9}{20} \cdot L(A,M) + \frac{11}{20} \cdot L(A,W) 
\hspace*{\arraycolsep} = \hspace*{\arraycolsep}
\frac{9}{20}\cdot \frac{1}{10} + \frac{11}{20}\cdot \frac{8}{10} 
\hspace*{\arraycolsep} = \hspace*{\arraycolsep}
\frac{97}{200}.
\end{array}$$

\noindent In the same way, $\subst{\graph(f)}(L)(B) = \frac{11}{20}$.  And
obviously,
$$\begin{array}{rcccccl}
\subst{\graph(f)}(L^\perp)(A)
& = &
\big(\subst{\graph(f)}(L)\big)^{\perp}(A)
& = &
1 - \subst{\graph(f)}(L)(A)
& = &
\frac{103}{200},
\end{array}$$

\noindent and similarly $\subst{\graph(f)}(L^\bot)(B) = \frac{9}{20}$.

\auxproof{
\begin{eqnarray*}
\subst{\graph(f)}(L)(B) & = & \sum_{(c,g) \in C \times G} \graph(f)(B)(c,g) \cdot L(c,g)\\
 & = & \frac{1}{2} \cdot L(B,M) + \frac{1}{2} \cdot L(B,W) \\
 & = & \frac{1}{2} \frac{2}{10} + \frac{1}{2}\frac{9}{10} \\
 & = & \frac{2}{20} + \frac{9}{20} \\
 & = & \frac{11}{20}
\end{eqnarray*}
}

With that done, we can now work out the conditional probabilities of a
man or a woman given the country and that they had long hair. In other
words, we are looking for a pair of maps $\given{f}{L},
\given{f}{L^\perp} \colon C \rightarrow G$ in $\Kl(\D)$ to fill in the
following triangle:
\begin{equation}
\label{exampleCondProbDiag}
\vcenter{\xymatrix@C+1pc{
& C\ar[dl]_{\charac_{\subst{\graph(f)}(L)} \quad}
   \ar[dr]^{\qquad j = (\pi_{2}+\pi_{2}) \klcomp \charac_{L} \klcomp \graph(f)} & \\
C+C\ar[rr]_-{(\given{f}{L}) + (\given{f}{L^{\perp}})} & & G+G
}}
\end{equation}

\noindent As distribution $\given{f}{L}(A) \in \D(G+G)$ we take:
$$\begin{array}{rcccccl}
\given{f}{L}(A)
& = &
{\displaystyle\frac{j(A) \after \kappa_{1}}{\subst{\graph(f)}(L)(A)}}
& = &
\frac{200}{97}\cdot \big(\frac{9}{200}\ket{M} + \frac{88}{200}\ket{W}\big)
& = &
\frac{9}{97}\ket{M} + \frac{88}{97}\ket{W}.
\end{array}$$

\noindent This can be read as: in country $A$, if we see someone with
long hair, the probability of this person being male (resp.\ female)
is $\frac{9}{97}$ (resp. $\frac{88}{97}$). In the same way one gets
$\given{f}{L}(B) = \frac{2}{11}\ket{M} + \frac{9}{11}\ket{W}$. And in
the two ``negated'' cases:
$$\begin{array}{rclcrcl}
\given{f}{L^\perp}(A)
& = &
\frac{81}{103}\ket{M} + \frac{22}{103}\ket{W} 
& \qquad &
\given{f}{L^\perp}(B)
& = &
\frac{8}{9}\ket{M} + \frac{1}{9}\ket{W}.
\end{array}$$

\auxproof{
\begin{eqnarray*}
f|L(A) & = & \frac{j(A)|_{\kappa_1G}}{\subst{\graph(f)}(L)(A)} \\
 & = & \frac{200}{97}\left(\frac{9}{200}\ket{M} + \frac{88}{200}\ket{W}\right) \\
 & = & \frac{9}{97}\ket{M} + \frac{88}{97}\ket{W} \\
f|L(B) & = & \frac{j(B)|_{\kappa_1G}}{\subst{\graph(f)}(L)(B)} \\
 & = & \frac{20}{11}\left(\frac{2}{20}\ket{M} + \frac{9}{20}\ket{W}\right) \\
 & = & \frac{2}{11}\ket{M} + \frac{9}{11}\ket{W}
\end{eqnarray*}

\begin{eqnarray*}
f|L^\bot(A) & = & \frac{j(A)|_{\kappa_2G}}{\subst{\graph(f)}(L^\bot)(A)} \\
 & = & \frac{200}{103} \left(\frac{81}{200}\ket{M} + \frac{22}{200}\ket{W} \right) \\
 & = & \frac{81}{103}\ket{M} + \frac{22}{103}\ket{W} \\
f|L^\bot(B) & = & \frac{j(B)|_{\kappa_2G}}{\subst{\graph(f)}(L^\bot)(B)} \\
 & = & \frac{20}{9}\left(\frac{8}{20}\ket{M} + \frac{1}{20}\ket{W} \right) \\
 & = & \frac{8}{9}\ket{M} + \frac{1}{9}\ket{W}
\end{eqnarray*}
}

\noindent With these definitions it is easy to see that
diagram~\eqref{exampleCondProbDiag} commutes.

\auxproof{
We now verify that these make diagram \ref{exampleCondProbDiag} commute.
\begin{eqnarray*}
& & (f|L + f|L^\bot)(\ch_{\subst{\graph(f)}(L)}(A)) \\
& = & (f|L + f|L^\bot)\left(\frac{97}{200}\ket{\kappa_1A} + \frac{103}{200}\ket{\kappa_2A}\right) \\
& = & \frac{97}{200}\kappa_1(f|L(A)) + \frac{103}{200}\kappa_2(f|L^\bot(A)) \\
& = & \frac{97}{200}\kappa_1\left(\frac{9}{97}\ket{M} + \frac{88}{97}\ket{W}\right) + \frac{103}{200}\kappa_2\left(\frac{81}{103}\ket{M} + \frac{22}{103}\ket{W} \right) \\
& = & \frac{9}{200}\ket{\kappa_1M} + \frac{88}{200}\ket{\kappa_1W} + \frac{81}{200}\ket{\kappa_2M} + \frac{22}{200}\ket{\kappa_2W} \\
& = & j(A)
\end{eqnarray*}

\begin{eqnarray*}
& & (f|L + f|L^\bot)(\ch_{\subst{\graph(f)}(L)}(B)) \\
& = & (f|L + f|L^\bot)\left(\frac{11}{20}\ket{\kappa_1B} + \frac{9}{20}\ket{\kappa_2B}\right) \\
& = & \frac{11}{20}\kappa_1(f|L(B)) + \frac{9}{20}\kappa_2(f|L^\bot(B)) \\
& = & \frac{11}{20}\kappa_1\left(\frac{2}{11}\ket{M} + \frac{9}{11}\ket{W}\right) + \frac{9}{20}\kappa_2\left(\frac{8}{9}\ket{M} + \frac{1}{9}\ket{W} \right) \\
& = & \frac{2}{20}\ket{\kappa_1M} + \frac{9}{20}\ket{\kappa_1W} + \frac{8}{20}\ket{\kappa_2M} + \frac{1}{20}\ket{\kappa_2W} \\
& = & j(B)
\end{eqnarray*}
}
\end{example}

Having completed an example, we move on to the general case for probabilities.

\begin{theorem}
\label{DiscProbThm}
For a morphism $f\colon X \rightarrow Y$ in the Kleisli category
$\Kl(\D)$ of the distribution monad $\D$, and for a predicate $\phi
\in [0,1]^{X\times Y}$, there are conditional probability maps
$\given{f}{\phi}, \given{f}{\phi^{\perp}} \colon X \rightarrow Y$ in
$\Kl(\D)$ making the following triangle commute.
$$\xymatrix@C+1pc{
& X\ar[dl]_{\charac_{\subst{\graph(f)}(\phi)} \quad }
   \ar[dr]^{\qquad \joint{f}{\phi} = (\pi_{2}+\pi_{2}) \klcomp \charac_{\phi} \klcomp \graph(f) } & \\
X+X\ar[rr]_-{(\given{f}{\phi}) + (\given{f}{\phi^{\perp}})} & & Y+Y
}$$

\noindent If $\subst{\graph(f)}(\phi)(x) \in (0,1)$, for each $x\in X$,
then both these maps $\given{f}{\phi}$ and $\given{f}{\phi^{\perp}}$
are uniquely determined.
\end{theorem}

\begin{myproof}
We define the functions $\given{f}{\phi}, \given{f}{\phi^{\perp}}
\colon X \rightarrow \D(Y)$ on $x\in X$ as:
$$\begin{array}{rclcrcl}
(\given{f}{\phi})(x)
& = &
\displaystyle\sum_{y} \frac{(\joint{f}{\phi})(x)(\kappa_{1}y)}
   {\subst{\graph(f)}(\phi)(x)}\ket{y}
& \qquad &
(\given{f}{\phi^{\perp}})(x)
& = &
\displaystyle\sum_{y} \frac{(\joint{f}{\phi})(x)(\kappa_{2}y)}
   {1-\subst{\graph(f)}(\phi)(x)}\ket{y}.
\end{array}$$

\noindent If $\subst{\graph(f)}(\phi)(x) \in \{0,1\}$, we choose an
arbitrary distribution instead. \QED
\end{myproof}

Here we have formulated conditional probability with respect to a
single formula $\phi$. It can be generalized to $n$-tests, which are
sequences of formulas $\phi_{1}, \ldots, \phi_{n}$ with $\phi_{1}
\ovee \cdots \ovee \phi_{n} = 1$ (see also~\cite{Jacobs13b}). Then one
gets $n$ corresponding conditional maps $\given{f}{\phi_i}$. In the
situation of the above theorem we actually use a 2-test, given by
$\phi$ and $\phi^{\perp}$.

\section{Conditional probability for $C^*$-algebras}\label{CstarSec}

We shall write $\CstarPU$ for the category of unital $C^*$-algebras (over the
complex numbers $\C$) with positive unital maps, and $\CstarMIU
\hookrightarrow \CstarPU$ for the subcategory where maps preserve
multiplication (M), involution (I), and unit (U); such maps are usually
called *-homomorphisms. In the present setting we assume all
$C^*$-algebras have a unit $1$. We write $\Ctr(A) \hookrightarrow A$
for the center of a $C^*$-algebra $A$, defined as usual as $\Ctr(A) =
\setin{a}{A}{\allin{b}{A}{ab = ba}}$. This center forms a commutative
(sub) $C^*$-algebra. Obviously, $A$ itself is commutative iff $A =
\Ctr(A)$. 

We remark at this point that any map in $\CstarPU$, when considered as
a map of Banach spaces, is of norm 1. This is \cite[corollary
  1]{RussoDye66}. This is equivalent to saying that any positive
unital map also preserves the norm.

The category $\CstarPU$ has finite products, via direct sums $\oplus$
of vector spaces (\textit{i.e.}~cartesian products of the underlying
sets). The operations are used pointwise. There are also tensor
products of $C^*$-algebras. These are described in more detail in
\cite[section IV.4]{TakesakiI}, but we outline them here. The
$C^*$-tensors for two $C^*$-algebras $A$ and $B$ are obtained by
taking the usual tensor of underlying vector spaces $A \otimes B$,
defining a *-algebra structure as follows:
$$\begin{array}{rclcrcl}
(a_1 \otimes b_1)(a_2 \otimes b_2) 
& = & 
(a_1a_2) \otimes (b_1b_2)
& \qquad &
(a \otimes b)^* 
& = &
a^* \otimes b^*.
\end{array}$$

\noindent One then obtains a $C^*$-algebra by introducing a $C^*$-norm
compatible with the *-algebra structure, and taking the
completion. There are minimal and maximal, or injective and projective
$C^*$-norms, but if $A$ or $B$ is finite dimensional these coincide
\cite[chapter XV, 1.4 - 1.6]{TakesakiIII}. Therefore, as a simplifying
assumption, all $C^*$-algebras in this paper will be taken as finite
dimensional unless explicitly stated otherwise.

The positive cone of $A \otimes B$ contains the positive elements
according to the multiplication and involution, \textit{i.e.}~$a \in A
\otimes B$ is positive if $a = b^*b$ for some other element $b$. We
note at this point that this cone is larger than the cone obtained by
taking sums of elements $a \otimes b$ with $a \in A$ and $b \in B$
both positive. The effect of this is that no $C^*$-tensor is a functor
on $\CstarPU$. The maps that can be tensored are called completely
positive and form a non-full subcategory $\CstarCPU \hookrightarrow
\CstarPU$ with the same objects (see \cite[section IV.3, and
  proposition IV.4.23]{TakesakiI}).

Since $A \otimes B$ is the completion of the algebraic tensor of $A$
and $B$, the span of elements of the form $a \otimes b$ is dense, and
in fact in our finite dimensional case $A \otimes B$ is just the span
of such elements. We can define coprojections $\kappa_i : A_i
\rightarrow A_1 \otimes A_2$ as follows:
\begin{align*}
\kappa_1(a) & = a \otimes 1 &
\kappa_2(a) & = 1 \otimes a,
\end{align*}

\noindent where $1$ is the unit of the $C^*$-algebra. It is simple to
see these are MIU maps, and therefore in $\CstarPU$.  

\auxproof{
\begin{eqnarray*}
\kappa_1(1) & = & 1 \otimes 1 \\
\kappa_1(a)\kappa_1(b) & = & (a \otimes 1)(b \otimes 1) \\
 & = & ab \otimes 1 \\
 & = & \kappa_1(ab) \\
\kappa_1(a)^* & = & (a \otimes 1)^* \\
 & = & a^* \otimes 1^* \\
 & = & a^* \otimes 1 \\
 & = & \kappa_1(a^*)
\end{eqnarray*}
}

(It is most natural to consider categories of $C^*$-algebras in
opposite form. For instance, in~\cite{FurberJ13a} it is shown that the
opposite $\op{(\CCstarPU)}$ of the category of commutative
$C^*$-algebras with positive unital maps is equivalent to a Kleisli
category, namely that of the ``Radon'' monad on compact Hausdorff
spaces. This restricts to an equivalence between finite-dimensional
commutative $C^*$-algebras and the subcategory $\KlN(\D)
\hookrightarrow \Kl(\D)$ with natural numbers as objects. In opposite
form, $\op{(\CstarCPU)}$ has similar structure to the Kleisli category
$\Kl(\D)$ used in the previous section, namely finite coproducts and
tensors with projections.)

We are working towards a $C^*$-algebraic analogue of
Lemma~\ref{KlDGraphLem}. But this requires some lemmas of its own. The
following result is based on theorem 1 of \cite{tomiyama}.
\begin{lemma}
\label{BimoduleLemma}
If $f : A \rightarrow B$ is a map in $\CstarMIU$.
\begin{enumerate}[(i)]
\item The algebra $B$ is a bimodule of $A$ under the left and right multiplications:
\begin{align*}
a \cdot b  & = f(a)b & b \cdot a & = b f(a).
\end{align*}
\item If a $\CstarPU$ map $g : B \rightarrow A$ is a retraction of
  $f$, \emph{i.e.} $g \after f = \id_A$, then $g$ is a map of
  bimodules:
\[
a_1g(b)a_2 = g(a_1 \cdot b \cdot a_2) = g(f(a_1)bf(a_2)).
\]
\end{enumerate}
\end{lemma}

\begin{myproof}
For the first point, the unit and multiplication properties follow
easily from those of $f$.  For the second point we notice that, since
$f$ has a left inverse, it is a split monic and therefore is
isomorphic to its image $f(A)$, a subalgebra of $B$. Then $f \after g$
is a positive unital projection onto $f(A)$, and is therefore a
projection of norm 1 in the sense of
\cite{tomiyama}. Applying~\cite[theorem 1, part 2]{tomiyama} we have
that $f\after g$ is a bimodule map. Thus, if $f(a_1), f(a_2) \in
f(A)$, and $b \in B$, then:
\[
(f\after g)\big(f(a_1)bf(a_2)\big) = f(a_1)(f \after g)(b)f(a_2) = f(a_1g(b)a_2),
\]

\noindent the latter because $f$ is a MIU-map. Applying the
injectivity of $f$, we have, as required:
$$\begin{array}{rcl}
g(f(a_1)bf(a_2))
& = &
a_1g(b)a_2.
\end{array}\eqno{\QEDbox}$$

\auxproof{
\begin{eqnarray*}
1 \cdot b & = & f(1)b \\
  & = & 1b = b \\
 & & \\
b \cdot 1 & = & bf(1) \\
 & = & b1 = b \\
 & & \\
a_1 \cdot (a_2 \cdot b) & = & f(a_1)(f(a_2)b) \\
 & = & f(a_1a_2)b \\
 & = & (a_1a_2)\cdot b
\end{eqnarray*}
The other proof is the mirror image of the above.
}

\end{myproof}

\begin{lemma}
\label{CstarGraphPosLem}
If $A$ is a $C^*$-algebra, multiplication of an element by an element
of the centre $\Ctr(A)$ is a MIU map $\mu\colon A \otimes \Ctr(A)
\rightarrow A$.
\end{lemma}

\begin{myproof}
Here is the definition of $\mu$:
$$\begin{array}{rcl}
\mu\left(\sum_ia_i \otimes z_i\right)
& = &
\sum_i a_iz_i.
\end{array}$$

\noindent Since the multiplication is bilinear, this is
well-defined. To show it preserves multiplication, it suffices to show it does so on basic tensors. We start with $\mu((a \otimes z)(b \otimes w)) = abzw$. Since $z$ commutes with $b$, we can rearrange this to get $azbw = \mu(a \otimes z) \mu(b \otimes w)$. The preservation of involution and unit are routine
arguments. \QED

\auxproof{ We now show that this is a MIU map:
\begin{itemize}
\item Multiplication: Let $a_1,a_2 \in A$ and $z_1,z_2 \in \Ctr(A)$.
\begin{eqnarray*}
\mu((a_1 \otimes z_1)\cdot(a_2 \otimes z_2)) & = & \mu(a_1a_2 \otimes z_1z_2) \\
 & = & a_1a_2z_1z_2
\end{eqnarray*}
Using the fact that $z_1$ commutes with $a_2$:
\begin{eqnarray*}
 & = & a_1z_1a_2z_2 \\
 & = & \mu(a_1\otimes z_1)\cdot \mu(a_2 \otimes z_2)
\end{eqnarray*}
\item Involution: Let $a \in A$, $z \in \Ctr(A)$.
\begin{eqnarray*}
\mu((a \otimes z)^*) & = & \mu(a^* \otimes z^*) \\
 & = & a^*z^*
\end{eqnarray*}
Using the fact that $z^*$ and $a^*$ commute:
\begin{eqnarray*}
 & = & z^*a^* \\
 & = & (az)^* \\
 & = & \mu(a \otimes z)^*
\end{eqnarray*}
\item Preservation of unit:
\begin{eqnarray*}
\mu(1 \otimes 1) & = & 1 \cdot 1 \\
 & = & 1
\end{eqnarray*}\qed
\end{itemize}
}
\end{myproof}

The following is the analogue of Lemma~\ref{KlDGraphLem}.

\begin{lemma}
\label{CstarGraphLem}
In $\CstarPU$ there is a bijective correspondence:
$$\begin{prooftree}
{\xymatrix{ B \ar[r]^-{f} & \Ctr(A) }}
\Justifies
{\xymatrix{ A\otimes B \ar[r]_-{g} & A \mbox{ positive and unital, with } 
   g \after \kappa_{1} = \idmap[A]}}
\end{prooftree}$$

\noindent Of course, when $A$ is commutative, the `$\Ctr$' can be dropped.

Like before we shall write $\graph(f)\colon A\otimes B \rightarrow A$
for the map corresponding to $f\colon B \rightarrow \Ctr(A)$, where $\graph(f) = \mu \after (\id_A \otimes f)$, or on elements
$\graph(f)(a\sotimes b) = a\cdot f(b) = f(b)\cdot a$. 
\end{lemma}

\begin{myproof}
Given $f : B \rightarrow \Ctr(A)$, in $\CstarPU$, since $\Ctr(A)$ is commutative we can use \cite[corollary IV.3.5]{TakesakiI} to show it is in $\CstarCPU$, and hence $\id_A \otimes f$ is positive. It is unital, and hence in $\CstarPU$ because $(\id_A \otimes f)(1 \otimes 1) = 1 \otimes f(1) = 1 \otimes 1$. By lemma \ref{CstarGraphPosLem}, $\mu$ is in $\CstarMIU$ and hence in $\CstarPU$, so $\mu \after (\id_A \otimes f)$ has the right type. To see it is a left inverse for $\kappa_1$:
\[
(\mu \after (\id_A \otimes f) \after \kappa_1)(a) 
 = \mu((\id_A \otimes f)(a \otimes 1))
 = \mu(a \otimes f(1))
 = \mu(a \otimes 1)
 = a.
\]
as required.

Conversely, if $g : A \otimes B \rightarrow A$ is a map such that $g \after \kappa_1 = \id_A$, then since $\kappa_1$ is an MIU map, Lemma \ref{BimoduleLemma} shows that $g$ is a bimodule map. We take $f$ to be $g \after \kappa_2 : B \rightarrow A$. This appears at first to have the wrong type. However, if $a \in A$ and $b \in B$
\[
a f(b) = a g(1 \otimes b) 
 = g(a\otimes 1 \cdot 1 \otimes b)
 = g(a \otimes b) 
 = g(1 \otimes b \cdot a \otimes 1) 
 = g(1 \otimes b) a 
 = f(b) a.
\]
Hence $f(b) \in \Ctr(A)$. It is left to the reader to check that the
correspondences we have described are each other's inverses. \QED

\auxproof{
$$\begin{array}{rcl}
\widehat{\graph(f)}(b)
& = &
\graph(f)(1\sotimes b) \\
& = &
1\cdot f(b) \\
& = &
f(b) \\
\graph(\widehat{g})(a\sotimes b)
& = &
a\cdot \widehat{g}(b) \\
& = &
a\cdot g(1\sotimes b) \\
& = &
g(a\sotimes b).
\end{array}$$
}
\end{myproof}

Also for $C^*$-algebras there is a logic $\op{(\CstarPU)} \rightarrow
\op{\EMod}$ of effect modules. For each $C^*$-algebra $A$, its
``effects'' $[0,1]_{A} = \setin{a}{A}{0 \leq a \leq 1}$ form an effect
module. The sum $e\ovee d$ exists and is equal to $e+d$ if $e+d\leq
1$.  The orthocomplement of $e$ is $e^{\perp} = 1 - e$. Each positive
unital map $f\colon A \rightarrow B$ forms an effect module map $\subst{f}
\colon [0,1]_{A} \rightarrow [0,1]_{B}$ by restriction. Since each
such map is determined by what it does on positive elements, we have a
full and faithful functor $\CstarPU \rightarrow \EMod$,
see~\cite{FurberJ13a} for more details.

The product $\times$ of $C^*$-algebras forms a coproduct $+$ in
$\op{(\CstarPU)}$. When we work in this opposite category, we shall
thus use the coproduct notation. There is a special effect $\Omega =
(1,0) \in [0,1]_{A+A} = ([0,1]_{A})^{2}$. For each effect $e\in
[0,1]_{A}$ there is a choice of characteristic map, in the non-commutative
case. In the category $\op{(\CstarPU)}$ one can define:
$$\xymatrix@C-.5pc{
A\ar[rr]^-{\charac_{e}} & & A+A
\qquad\mbox{as}\qquad
(a,a')\ar@{|->}[r] & 
    \sqrt{e}\cdot a \cdot \sqrt{e} + \sqrt{(1-e)}\cdot a'\cdot \sqrt{(1-e)} \\
}$$

\noindent We have that $\subst{\charac_{e}}(\Omega) =
\charac_{e}(\Omega) = e$, and $\subst{\charac_{e}}(\Omega^\bot) = \charac_e(\Omega^\bot) = e^\bot$. This property replaces that of \cite{newdirections}, being a section of $\nabla_A$, which is not satisfied in the non-commutative case.

\begin{lemma}
\label{TwoSideProductLemma}
If $a,b$ are positive elements of a $C^*$-algebra, then $aba$ is positive.
\end{lemma}
\begin{myproof}
Since $b$ is positive, $b = p^*p$ for some $p$. So
\[
aba = ap^*pa = a^*p^*pa = (pa)^*pa
\]
and $aba$ is positive. (We have used that all positive elements are self-adjoint.) \QED
\end{myproof}

\begin{corollary}
\label{CharMapPosCor}
The map $\charac_e$ is positive and unital.
\end{corollary}

\begin{myproof}
Effects $e$ and $e^{\perp} = 1-e$ are positive, and so are their
square roots.  Hence the previous lemma makes $\charac_e$
positive. The proof of unitality is straightforward. \QED

\auxproof{
\begin{itemize}
\item Unitality:
\begin{eqnarray*}
\charac_e(1,1) & = & \sqrt{e}\cdot 1 \cdot \sqrt{e} + \sqrt{(1-e)}\cdot 1 \cdot \sqrt{(1-e)} \\
 & = & e + 1 - e = 1
\end{eqnarray*}
\item Positivity:
Since $\sqrt{e}$ is defined to be the positive square root, by lemma \ref{TwoSideProductLemma}
\[
\charac_e(a_1,a_2) = \sqrt{e}a_1\sqrt{e} + \sqrt{(1-e)}a_2\sqrt{(1-e)}
\]
is the sum of two positive elements and therefore positive. \QED
\end{itemize}
}
\end{myproof}

We can now give a proof that this definition of characteristic map, for commutative $C^*$-algebras, coincides with the monadic definition for the Radon monad under the equivalence $\Kl(\Rdn) \simeq \CCstarPU$ from \cite[Theorem 2]{FurberJ13a}. In both cases we can start with a compact Hausdorff space $X$, and take the corresponding $C^*$-algebra to be $C(X)$, the $C^*$-algebra of continuous functions $X \rightarrow \C$. A predicate is a continuous map to the unit interval, $e \in \CHaus(X, [0,1])$. We have two possible characteristic maps
\[
\ch_e : C(X) \times C(X) \rightarrow C(X)
\]
and
\[
\ch_e' : X \rightarrow \Rdn(X + X)
\]
which, following \eqref{CharacProbEqn}, is defined as
\[
\ch_e'(x) = e(x)\delta_{\kappa_1 x} + (1-e(x))\delta_{\kappa_2 x},
\]
where the $\delta$s are Dirac delta measures. This may equivalently be defined, given a function $f \in C(X + X)$, as
\[
\ch_e'(x)(f) = e(x) \cdot f(\kappa_1 x) + (1-e(x)) \cdot f(\kappa_2 x) .
\]

\begin{theorem}
Under the equivalence $\CR : \Kl(\Rdn) \tilde{\rightarrow} \CCstarPU$ from \cite[Theorem 2]{FurberJ13a}, $\ch_e$ coincides with $\ch_e'$, which is to say, given $a_1, a_2 \in C(X)$
\[
\CR(\ch_e')([a_1,a_2]) = \ch_e(a_1,a_2).
\]
\end{theorem}
\begin{proof}
Consider the right hand side. We have that 
\[
\charac_e(a_1,a_2) = \sqrt{e}a_1\sqrt{e} + \sqrt{1-e}a_2\sqrt{1-e},
\]
which by commutativity can be rewritten as
\[
\charac_e(a_1,a_2) = ea_1 + (1-e)a_2 .
\]
Now let $x \in X$, and we can see
\begin{align*}
\CR(\ch_e')([a_1,a_2])(x) &= \ch_e'(x)([a_1,a_2]) \\
 &= e(x) \cdot [a_1,a_2](\kappa_1 x) + (1-e(x)) \cdot [a_1,a_2](\kappa_2 x) \\
 &= e(x) \cdot a_1(x) + (1-e(x)) \cdot a_2(x) \\
 &= (e a_1 + (1-e) a_2)(x)
\end{align*}
and so $\CR(\ch_e')([a_1,a_2]) = \ch_e(a_1,a_2)$ as required.
\end{proof}

\auxproof{
The following is a proof that the two definitions of characteristic map agree on commutative C$^*$-algebras, using the Radon monad analogue of the first one.

First we rewrite the characteristic maps as follows:

In the C$^*$-case
\[
\charac_e(a_1,a_2) = \sqrt{e}a_1\sqrt{e} + \sqrt{1-e}a_2\sqrt{1-e}
\]
which by commutatitivity can be rewritten as
\[
\charac_e(a_1,a_2) = ea_1 + (1-e)a_2
\]

and for $\Kl(\Rdn)$, we have
\[
\charac_e : X \rightarrow \Rdn(X + X)
\]

which is defined as
\[
\charac_e(x)(f) = e(x)f(\kappa_1(x)) + (1-e(x))\kappa_2(x)
\]

If we consider $S(\charac_e) : S(A) \rightarrow S(A \times A)$ we can get the Kleisli map by evaluating at delta functions.

\begin{align*}
S(\charac_e)(\delta_x)(f : X + X \rightarrow \C) &= (\delta_x \circ \charac_e)(f) \\
 &= \delta_x(\charac_e(f)) \\
 &= \delta_x(\charac_e(f \circ \kappa_1, f \circ \kappa_2)) \\
 &= \delta_x(e \cdot (f \circ \kappa_1) + (1-e)\cdot (f \circ \kappa_2)) \\
 &= (e \cdot (f \circ \kappa_1) + (1-e) \cdot (f \circ \kappa_2))(x) \\
 &= e(x)\cdot f(\kappa_1(x)) + (1 - e(x))\cdot f(\kappa_2(x)) \\
 &= \charac_e(x)(f)
\end{align*}

}

We are now in a position to describe a setting for conditional
probability for $C^*$-algebras. In order to maximize the analogy with
the situation in the previous section --- involving the Kleisli
category $\Kl(\D)$ --- we work in the opposite category
$\op{(\CstarPU)}$. There tensors have projections $\pi_i$.

Assume we have a map $f\colon \Ctr(A) \rightarrow B$ and an effect
$e\in [0,1]_{A\otimes B}$. Then we can form the marginal and total
probability maps as follows.
\begin{itemize}
\item Via the graph map $\graph(f) \colon A \rightarrow A\otimes B$
obtained in Lemma~\ref{CstarGraphLem} we can substitute and get
$\subst{\graph(f)}(e) \in [0,1]_{A}$ and form the characteristic map
$\charac_{\subst{\graph(f)}(e)} \colon A \rightarrow A+A$.

\item We can also form the joint probability $\joint{f}{e}$ as the
  composition, in $\op{(\CstarPU)}$:
$$\xymatrix@C+.5pc{
\joint{f}{e} = \Big(A\ar[r]^-{\graph(f)} & 
   A\otimes B\ar[r]^-{\charac_{e}} & 
   (A\otimes B)+(A\otimes B)\ar[r]^-{\pi_{2}+\pi_{2}} & B+B\Big)
}$$
\end{itemize}

\noindent The conditional probability maps $\given{f}{e}, \given{f}{e^{\perp}}
\colon A \rightarrow B$ in $\op{(\CstarPU)}$ then fit in the triangle:
\begin{equation}
\label{QuantumParamCondExpDiag}
\vcenter{\xymatrix@C+1pc{
& A\ar[dl]_{\charac_{\subst{\graph(f)}(e)} \quad }
   \ar[dr]^{\quad \joint{f}{e}} & \\
A+A\ar[rr]_-{(\given{f}{e}) + (\given{f}{e^{\perp}})} & & B+B
}}
\end{equation}

\begin{theorem}
\label{ParamCondExpTheorem}
If $\subst{\graph(f)}(e)$ and $\subst{\graph(f)}(e^\bot)$ are invertible, then the maps $\given{f}{e}$ and $\given{f}{e^\bot}$ exist and are unique. The formulas for each are:
\begin{eqnarray*}
\given{f}{e}(b) & = & \isqrt{\subst{\graph(f)}(e)} \cdot 
   (\joint{f}{e})(b,0) \cdot \isqrt{\subst{\graph(f)}(e)} \\
\given{f}{e^\bot}(b) & = & \isqrt{\subst{\graph(f)}(1-e)} \cdot 
   (\joint{f}{e})(0,b) \cdot \isqrt{\subst{\graph(f)}(1-e)} 
\end{eqnarray*}
\end{theorem}

\begin{myproof}
The proof has three steps. First we show that these maps are in
$\CstarPU$. Then we show they make the diagram commute. Finally, we
show they are the unique such maps.

But first, we remark that for any positive invertible element $a$ of a
$C^*$-algebra $A$, the spectrum of $a$ is a closed subset of
$(0,\infty)$ and so we may use continuous functional calculus (see
\cite[definition I.4.7]{TakesakiI}) to take $\isqrt{a}$, so the
positive square root of $a$ is invertible.

These maps can be seen to be unital because the inverse square roots
on either side cancel with the square roots.  

\auxproof{
\begin{eqnarray*}
\given{f}{e}(1) & = & \isqrt{\subst{\graph(f)}(e)} \cdot 
   (\joint{f}{e})(1,0) \cdot \isqrt{\subst{\graph(f)}(e)} \\
 & = & \isqrt{\subst{\graph(f)}(e)} \cdot \subst{\graph(f)}(\charac_e(\kappa_2 \times \kappa_2(1,0))) \cdot \isqrt{\subst{\graph(f)}(e)} \\
 & = & \isqrt{\subst{\graph(f)}(e)} \cdot \subst{\graph(f)}(\charac_e(1 \otimes 1,0)) \cdot \isqrt{\subst{\graph(f)}(e)} \\
 & = & \isqrt{\subst{\graph(f)}(e)} \cdot \subst{\graph(f)}(e) \cdot \isqrt{\subst{\graph(f)}(e)} \\
 & = & \isqrt{\subst{\graph(f)}(e)} \cdot \sqrt{\subst{\graph(f)}(e)}\cdot \sqrt{\subst{\graph(f)}(e)} \cdot \isqrt{\graph(\subst{f}(e))}\\
 & = & 1
\end{eqnarray*}
The other map's proof is similar.
}

To show that they are positive, let $a$ be a positive element of $A$. Then:
$$\begin{array}{rcl}
\given{f}{e}(a)
& = &
\isqrt{\subst{\graph(f)}(e)} \cdot 
  \subst{\graph(f)}\big(\charac_e\big((\kappa_2 \times \kappa_2)(a,0)\big)\big)    \cdot \isqrt{\subst{\graph(f)}(e)}.
\end{array}$$

\noindent If we show that $(\kappa_2 \times \kappa_2)(a,0)$ is
positive, then it will follow from corollary~\ref{CharMapPosCor},
lemma~\ref{CstarGraphLem} and lemma~\ref{TwoSideProductLemma} that
$\given{f}{e}(a)$ is positive. Since $(\kappa_2 \times \kappa_2)(a, 0)
= (1 \otimes a, 0)$ and $a$ is positive, it can be written $a = b^*b$,
so that we have:
\[
(1 \otimes b^*b, 0) = ((1 \otimes b^*)(1 \otimes b), 0) = ((1 \otimes b)^*(1 \otimes b), 0),
\]
which is positive. Thus $\given{f}{e}(a)$ is positive. The case of $\given{f}{e^\bot}$ is similar.

To show that these maps $\given{f}{e}$, $\given{f}{e^\perp}$ make the
diagram~\eqref{QuantumParamCondExpDiag} commute, let $(b_1,b_2) \in B
\times B$, where we are reconsidering the diagram in $\CstarPU$. Then
$$\begin{array}{rcl}
\lefteqn{\charac_{\subst{\graph(f)}(e)}\big(
   (\given{f}{e} \times \given{f}{e^\bot})(b_1,b_2)\big)} \\
& = & 
\charac_{\subst{\graph(f)}(e)}\big((\given{f}{e})(b_1), (\given{f}{e^\bot})(b_2)\big) \\
& = & 
\charac_{\subst{\graph(f)}(e)} \left( \isqrt{\subst{\graph(f)}(e)} \cdot 
   (\joint{f}{e})(b_1,0) \cdot \isqrt{\subst{\graph(f)}(e)}, 
   \isqrt{\subst{\graph(f)}(1-e)} \cdot 
   (\joint{f}{e})(0,b_2) \cdot \isqrt{\subst{\graph(f)}(1-e)} \right) \\
& = & 
(\joint{f}{e})(b_1,0) + (\joint{f}{e})(0,b_2) \\
& = & 
(\joint{f}{e})((b_1,0) + (0,b_2)) \\
 & = & 
(\joint{f}{e})(b_1,b_2)
\end{array}$$

To show the uniqueness, suppose we have $g_1,g_2$ such that
$\charac_{\subst{\graph(f)}(e)} \after (g_1 \times g_2) =
\joint{f}{e}$, and let $(b_1,b_2) \in B \times B$.  Then if $b_1 \in
B$ we have:
$$\begin{array}{rcl}
\charac_{\subst{\graph(f)}(e)}\big((g_1 \times g_2)(b_1,0)\big)
& = &
(\joint{f}{e})(b_1,0).
\end{array}$$

\noindent Rearranging the left hand side, we get:
$$\begin{array}{rcccl}
\charac_{\subst{\graph(f)}(e)}(g_1(b_1),0) 
& = & 
\sqrt{\subst{\graph(f)}(e)} \cdot g_1(b_1) \cdot 
   \sqrt{\subst{\graph(f)}(e)} + 0 
 & = & 
(\joint{f}{e})(b_1,0).
\end{array}$$

\noindent By invertibility of $\sqrt{\subst{\graph(f)}(e)}$, we have
that:
$$\begin{array}{rcl}
g_1(b_1)
& = &
\isqrt{\subst{\graph(f)}(e)} \cdot (\joint{f}{e})(b_1,0) \cdot 
   \isqrt{\subst{\graph(f)}(e)},
\end{array}$$

\noindent as required. The $g_2$ case is similar. \QED
\end{myproof}

For ease of application in the example in the next section, we
specialize quantum conditional probability to when there is no
parametrization, taking $A = \C$ in
Theorem~\ref{ParamCondExpTheorem}. Then $f$ is a state, considered as
a map $B \rightarrow \C$ in $\CstarPU$. We can use the isomorphism of
$B \otimes \C \cong B$ to view $e$ as a predicate on $B$. Then
$\joint{f}{e} = f \after \charac_e$, much like
in~\eqref{HairTriangle}. This means
diagram~\eqref{QuantumParamCondExpDiag} becomes in $\op{(\CstarPU)}$:
\begin{equation}
\label{QuantumCondExpDiag}
\vcenter{\xymatrix@R+1pc@C+1pc{
& \C\ar[dl]_{\charac_{\subst{f}(e)}}
   \ar[dr]^{\charac_{e} \after f} & \\
\C+\C\ar[rr]_-{(\given{f}{e})+ (\given{f}{e^{\perp}})} & & B + B
}}
\end{equation}

\begin{corollary}
For a state $f\colon B \rightarrow \C$ and an effect $e\in [0,1]_{B}$,
if $\subst{f}(e) = f(e) \neq 0,1$ then the conditional states
$\given{f}{e}, \given{f}{e^\perp}$ in~\eqref{QuantumCondExpDiag} exist
and are unique, and can be given by the formulas:
\begin{equation}
\label{UnParamGivenEqn}
\begin{array}{rclcrcl}
\given{f}{e}(b) 
& = & 
\displaystyle\frac{f\big(\sqrt{e}b\sqrt{e}\big)}{f(e)}
& \qquad &
\given{f}{e^\bot}(b) 
& = & 
\displaystyle\frac{f\big(\sqrt{1-e}b\sqrt{1-e}\big)}{f(1-e)}.
\end{array}
\end{equation}
\end{corollary}

\begin{myproof}
Since $\C$ is a field, $0$ is the only non-invertible element. Since
$\subst{f}(1-e) = 1 - \subst{f}(e)$ as $f$ is unit-preserving and
linear, $\subst{f}(e) \neq 0,1$ implies that $\subst{f}(e)$ and
$\subst{f}(1-e)$ are invertible. We then apply
Theorem~\ref{ParamCondExpTheorem} and use the commutativity of
$\C$. \QED

\auxproof{ Theorem
  \ref{ParamCondExpTheorem} produces:
\begin{eqnarray*}
\given{f}{e}(b) & = & \isqrt{\subst{f}(e)} \cdot \joint{f}{e}(b,0) \cdot \isqrt{\subst{f}(e)} \\
\given{f}{e^\bot}(b) & = & \isqrt{\subst{f}(1-e)} \cdot \joint{f}{e}(0,b) \cdot \isqrt{\subst{f}(1-e)}
\end{eqnarray*}
which by the commutativity of $\C$ implies
\begin{eqnarray*}
\given{f}{e}(b) & = & \frac{\joint{f}{e}(b,0)}{\subst{f}(e)} \\
\given{f}{e^\bot}(b) & = & \frac{\joint{f}{e}(0,b)}{\subst{f}(1-e)}
\end{eqnarray*}
Using the equivalent form of the joint probability in this case gives
\begin{eqnarray*}
\given{f}{e}(b) & = & \frac{f(\ch_e(b,0))}{\subst{f}(e)} \\
\given{f}{e^\bot}(b) & = & \frac{f(\ch_e(0,b))}{\subst{f}(1-e)}
\end{eqnarray*}
as is required.}
\end{myproof}

Since this definition of conditional probability applies to effects, not just projections, it in fact works as a definition of conditional expectation for postiive operators less than or equal to $1$. As such, it may be related to the definition of conditional expectation given in \cite{chaput}. However, we have used $C^*$-algebras here and a non-commutative version of the definition in that paper is more naturally formulated in the setting of $W^*$-algebras, so we leave relating the two to future work.

\section{Example}

As an example, we use the bomb tester of \cite{Elitzur93}. Suppose
some bombs exist that explode if a single photon is absorbed by a
detector attached to it. However, some of these bombs are duds, and
the photon passes through the detector unaltered, failing to explode
the bomb, if this is the case. We want to find out which of the bombs
are which. If we try to test a bomb to see if it explodes, we
seemingly can only keep the bomb if it turns out to be a dud, as the
bomb will explode if tested with a photon, the smallest amount of
light that we could use. However it is shown in~\cite{Elitzur93} that
this is not the case, and a bomb tester can be built. We reformulate
this to use our framework for quantum conditional probability.

\begin{figure}[htbp]
\label{bombfig}
\begin{center}
\scalebox{0.5}{
\includegraphics{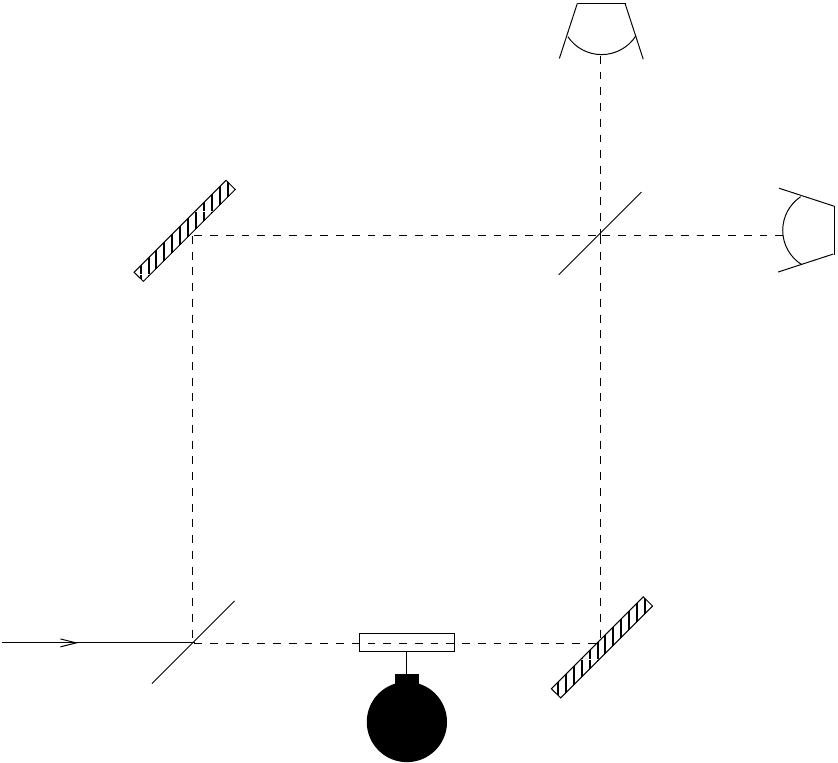}
}
\caption{The bomb tester}
\end{center}
\end{figure}

The set-up is similar to a Mach-Zehnder interferometer, as observed
in~\cite{Elitzur93}. A photon passes through a semi-silvered mirror,
where the bomb is in the path of one branch, the photon is reflected
from two mirrors to hit a second semi-silvered mirror, after which
there are two detectors. This can be seen in figure~\ref{bombfig}. We
represent the system with the following $C^*$-algebra:
$$\begin{array}{rcccl}
A
& = &
A_{E} \otimes A_{P} \otimes A_{B} 
& = &
C(\{\mathrm{L},\mathrm{D}\}) \otimes 
   B(\ell^2(\{\uparrow, \rightarrow, \emptyset\})) \otimes B(\ell^2(\{0,1\})).
\end{array}$$

\noindent The status of the bomb being Live or a Dud is treated as
classical, the direction or absence of a photon is represented by a
3-dimensional Hilbert space and the state of the bomb as unexploded or
exploded is treated as a 2-dimensional Hilbert space. We use the
shortened names $A_E, A_P$ and $A_B$ for these algebras, the letters
standing for explosivity, photon, and bomb respectively. All together,
the $C^*$-algebra is $2 \times 3^2 \times 2^2 = 72$-dimensional.

The mirrors (semi-silvered or fully silvered) act only on $A_P$. They
are maps of the form $U^*\cdot - \cdot U$ for $U$ a unitary from
$\ell^2(\{\uparrow, \rightarrow, \emptyset\})$ to itself.  On basis
vectors, the semi-silvered mirrors' unitaries, $U_S$, are:
$$\begin{array}{rclcrclcrcl}
\ket{\rightarrow} 
& \mapsto & 
\frac{1}{\sqrt{2}}\ket{\rightarrow} + \frac{1}{\sqrt{2}}\ket{\uparrow}
& \qquad &
\ket{\uparrow} 
& \mapsto & 
\frac{1}{\sqrt{2}}\ket{\rightarrow} - \frac{1}{\sqrt{2}}\ket{\uparrow}
& \qquad &
\ket{\emptyset} 
& \mapsto & 
\ket{\emptyset}.
\end{array}$$

\noindent And the fully silvered mirrors' unitaries, $U_F$, are:
$$\begin{array}{rclcrclcrcl}
\ket{\rightarrow} & \mapsto & \ket{\uparrow}
& \qquad &
\ket{\uparrow} & \mapsto & \ket{\rightarrow}
& \qquad &
\ket{\emptyset} & \mapsto & \ket{\emptyset}.
\end{array}$$

\noindent The reader may verify that these are unitary and that
$U_SU_FU_S\ket{\rightarrow} = \ket{\rightarrow}$, so that in the
absence of a bomb the photon always comes out to the right. Already, a
stark difference is apparent from what would happen if the
semi-silvered mirrors acted probabilistically.

The explosion of the bomb can be represented as a unitary $U_B$ in
$A_P \otimes A_B$. To do this, we allow the bomb to spontaneously
unexplode, emitting a rightward photon. This does not affect the
results as the map is never evaluated in this state.
$$\begin{array}{rclcrclcrcll}
\ket{\uparrow 0} & \mapsto & \ket{\uparrow 0} 
& \qquad & 
\ket{\uparrow 1} & \mapsto & \ket{\uparrow 1} 
& \qquad &
\ket{\rightarrow 0} & \mapsto & \ket{\emptyset 1} 
& \mbox{(bomb explodes)} \\
\ket{\rightarrow 1} & \mapsto & \ket{\rightarrow 1}
& &
\ket{\emptyset 0} & \mapsto & \ket{\emptyset 0}
& &
\ket{\emptyset 1} & \mapsto & \ket{\rightarrow 0} 
   &\mbox{(bomb unexplodes).}
\end{array}$$

\noindent We can then describe the dynamics of the exploding bomb on
$A = A_E \otimes A_P \otimes A_B$. For ease of use later, we actually
use the Schr\"odinger picture. We have the states $\delta_x \in
\CstarPU(A_{E}, \C) \cong \D(\{L, D\})$, for $x \in \{\mathrm{L},
\mathrm{D}\}$. It is the usual delta measure, which is a map $A_E
\rightarrow \C$. For a state $\rho : A_P \otimes A_B \rightarrow \C$
the $\delta_x$ determine the dynamics, as in:
$$\begin{array}{rclcrcl}
\delta_{\mathrm{L}} \otimes \rho 
& \mapsto & 
\delta_{\mathrm{L}} \otimes \rho(U_B^* \cdot \mbox{ - } \cdot U_B)
& \qquad &
\delta_{\mathrm{D}} \otimes \rho 
& \mapsto & 
\delta_{\mathrm{D}} \otimes \rho.
\end{array}$$

\noindent As we can see, whether the bomb can explode or not depends
on whether we have an $L$ or $D$ state in the first component of $A$.

The way conditional probability is supposed to work is that for $x \in
A$ considered to be a random variable, we have some $f : A \rightarrow
\C$ such that $f(x) = \E(x)$, and $\given{f}{e}(x) = \E(x \mid e =
1)$, as in diagram \eqref{QuantumCondExpDiag}. To work out $f$, we
start off with the initial state:
$$\xymatrix@C+2pc{
f_0 = \Big(A_{E}\otimes (A_{P}\otimes A_{B})
   \ar[rr]^-{(\frac{1}{2}\delta_{\mathrm{D}} + \frac{1}{2}\delta_{\mathrm{L}})
   \otimes \bra{\rightarrow 0} \mbox{ - } \ket{\rightarrow 0}} & &
   \C\otimes \C \cong \C\Big).
}$$

\noindent In other words, we start with an even probability of a live
bomb or a dud, and with the photon moving to the right, before it hits
the first mirror.

Then we show how $f_0$ changes under the dynamics. To save space, we replace 
the ket of the state with an ellipsis (\ldots):
\begin{eqnarray*}
\mbox{first mirror} & & (\frac{1}{2}\delta_{\mathrm{D}} + \frac{1}{2}\delta_{\mathrm{L}})\otimes (\frac{1}{\sqrt{2}}\bra{\rightarrow 0} + \frac{1}{\sqrt{2}}\bra{\uparrow 0}) \mbox{ - } \ket{\cdots} \\
\mbox{light hits bomb} & & \frac{1}{2}\delta_{\mathrm{D}}\otimes(\frac{1}{\sqrt{2}}\bra{\rightarrow 0} + \frac{1}{\sqrt{2}}\bra{\uparrow 0})\mbox{ - }\ket{\cdots} + \frac{1}{2}\delta_{\mathrm{L}} \otimes (\frac{1}{\sqrt{2}}\bra{\emptyset 1} + \frac{1}{\sqrt{2}}\bra{\uparrow 0}) \mbox{ - } \ket{\cdots} \\
\mbox{opaque mirrors} & & \frac{1}{2}\delta_{\mathrm{D}}\otimes(\frac{1}{\sqrt{2}}\bra{\uparrow 0} + \frac{1}{\sqrt{2}}\bra{\rightarrow 0})\mbox{ - }\ket{\cdots} + \frac{1}{2}\delta_{\mathrm{L}} \otimes (\frac{1}{\sqrt{2}}\bra{\emptyset 1} + \frac{1}{\sqrt{2}}\bra{\rightarrow 0}) \mbox{ - } \ket{\cdots} \\
\mbox{last mirror} & & \frac{1}{2}\delta_{\mathrm{D}}\otimes \bra{\rightarrow 0}\mbox{ - }\ket{\rightarrow 0} + \frac{1}{2}\delta_{\mathrm{L}} \otimes (\frac{1}{\sqrt{2}}\bra{\emptyset 1} + \frac{1}{2}\bra{\rightarrow 0} + \frac{1}{2}\bra{\uparrow 0}) \mbox{ - } \ket{\cdots} 
\end{eqnarray*}

\noindent We shall write $f$ for this last state $A
\rightarrow \C$.  Now that it is fixed, consider the situation in which
the bomb did not explode and the photon was detected going up. This
is captured by the following effect.
$$\begin{array}{rcl}
e & = & 1_{A_E} \otimes \ketbraI{\uparrow 0} \;\in\; [0,1]_{A}.
\end{array}$$

\noindent We wish to calculate the probability that the bomb is a dud,
given $e$, \textit{i.e.}~given that the bomb did not explode and the
photon was detected going up. In symbols this is $\Prob(\mbox{Dud}
\mid e)$, \textit{i.e.} $\E(\chi_D \otimes 1_{A_P \otimes A_B} \mid
e)$, where $\chi_{D}\in C(\{L,D\})$ is the obvious indicator
function. In triangle diagram~\eqref{QuantumCondExpDiag} we wish to
calculate the conditional state $\given{f}{e}$ with input event $b =
\chi_D \otimes 1_{A_P \otimes A_B} \in A$.

We apply the formula~\eqref{UnParamGivenEqn} for $\given{f}{e}(b) \in
\C$. To do this, we first calculate $\subst{f}(e) = f(e)$, using the
abbreviation $\ket{\psi} = \frac{1}{\sqrt{2}}\ket{\emptyset 1} +
\frac{1}{2} \ket{\rightarrow 0} + \frac{1}{2} \ket{\uparrow 0}$. First
$\subst{f}(e)$:
$$\begin{array}{rcl}
f(e) 
\hspace*{\arraycolsep} = \hspace*{\arraycolsep}
f(1_{A_E} \otimes \ketbraI{\uparrow 0}) 
& = & 
\frac{1}{2}\delta_{\mathrm{D}}(1)\otimes \braket{\rightarrow 0}{\uparrow 0} 
  \braket{\uparrow 0}{\rightarrow 0} + \frac{1}{2}\delta_{\mathrm{L}}(1) \otimes 
  \braket{\psi}{\uparrow 0}\braket{\uparrow 0}{\psi} \\
& = & 
0 + \frac{1}{2} (\frac{1}{2}\braketI{\uparrow 0} 
   \braketI{\uparrow 0}\frac{1}{2}) 
\hspace*{\arraycolsep} = \hspace*{\arraycolsep}
\frac{1}{8}.
\end{array}$$

\noindent Since the effect $e$ is a projection, it is its own positive
square root. Therefore we have:
\[
\charac_e(b,0) 
=
\sqrt{e}b\sqrt{e}
=
ebe
=
(1_{A_E} \otimes \ketbraI{\uparrow 0}) (\chi_{\mathrm{D}} \otimes 1_{A_P \otimes A_B}) 
   (1_{A_E} \otimes \ketbraI{\uparrow 0}) 
=
\chi_{\mathrm{D}} \otimes \ketbraI{\uparrow 0}.
\]

\noindent We may now substitute all of these values
into~\eqref{UnParamGivenEqn} and get $\given{f}{e}(b)$:
$$\begin{array}{rcccccl}
\given{f}{e}(b) 
& = & 
\displaystyle\frac{f(\chi_{\mathrm{D}} \otimes \ketbraI{\uparrow 0})}
   {\frac{1}{8}}
 & = & 
8\Big(\frac{1}{2}\braket{\rightarrow 0}{\uparrow 0} 
   \braket{\uparrow 0}{\rightarrow 0} + 0\Big)
 & = & 
0.
\end{array}$$

Thus if an upward-moving photon is detected and the bomb did not explode, the probability that it is a dud is $0$, and it must be live. This gives a way to get live bombs without exploding them. Note that this contradicts a commonly stated notion about quantum mechanics, that one cannot observe something without affecting it\footnote{This is intended to refer to the projection that occurs in a measurement}, as in this case we have a way to use quantum mechanics to observe something without affecting it in a way that we would have had to do classically.

\section{Conclusions}

In this paper we have given a categorical formulation of conditional
probability. It involves a triangle-fill property, where the condition
is a predicate from an associated predicate logic, formalized via an
indexed category of effect modules. It is shown that this formulation
gives the familiar classical notion of conditional probability, when
interpreted in the Kleisli category of the distribution monad.

Next, the formulation can also be used in a quantum setting, given by
the category of finite-dimensional $C^*$-algebras. We have presented a
general ``parametrized'' formulation, but our main example, the bomb
tester, only involves the non-parametrized case. Further clarification
is needed, in this general parametrized case, also in relation to
other approaches in the literature. Our approach has the advantage
that it is based on a general categorical scheme, that can be
instantiated in various settings.

\subsubsection*{Acknowledgements}
This research has been financially supported by the Netherlands Organisation for Scientific Research (NWO) under TOP-GO grant no. 613.001.013 (The logic of composite quantum systems).

\bibliographystyle{eptcs}
\bibliography{bayes}

\begin{thebibliography}{10}
\providecommand{\bibitemdeclare}[2]{}
\providecommand{\surnamestart}{}
\providecommand{\surnameend}{}
\providecommand{\urlprefix}{Available at }
\providecommand{\url}[1]{\texttt{#1}}
\providecommand{\href}[2]{\texttt{#2}}
\providecommand{\urlalt}[2]{\href{#1}{#2}}
\providecommand{\doi}[1]{doi:\urlalt{http://dx.doi.org/#1}{#1}}
\providecommand{\bibinfo}[2]{#2}

\bibitemdeclare{article}{bub07}
\bibitem{bub07}
\bibinfo{author}{J.~\surnamestart Bub\surnameend} (\bibinfo{year}{2007}):
  \emph{\bibinfo{title}{{Quantum Probabilities as Degrees of Belief}}}.
\newblock {\sl \bibinfo{journal}{Studies in History and Philosophy of Science
  Part B: Studies in History and Philosophy of Modern Physics}}
  \bibinfo{volume}{38}(\bibinfo{number}{2}), pp. \bibinfo{pages}{232 -- 254},
  \doi{10.1016/j.shpsb.2006.09.002}.

\bibitemdeclare{incollection}{chaput}
\bibitem{chaput}
\bibinfo{author}{P.~\surnamestart Chaput\surnameend},
  \bibinfo{author}{V.~\surnamestart Danos\surnameend},
  \bibinfo{author}{P.~\surnamestart Panangaden\surnameend} \&
  \bibinfo{author}{G.~\surnamestart Plotkin\surnameend} (\bibinfo{year}{2009}):
  \emph{\bibinfo{title}{Approximating Markov Processes by Averaging}}.
\newblock In \bibinfo{editor}{S.~\surnamestart Albers\surnameend},
  \bibinfo{editor}{A.~\surnamestart Marchetti-Spaccamela\surnameend},
  \bibinfo{editor}{Y.~\surnamestart Matias\surnameend},
  \bibinfo{editor}{S.~\surnamestart Nikoletseas\surnameend} \&
  \bibinfo{editor}{W.~\surnamestart Thomas\surnameend}, editors: {\sl
  \bibinfo{booktitle}{Automata, Languages and Programming}}, {\sl
  \bibinfo{series}{Lecture Notes in Computer Science}} \bibinfo{volume}{5556},
  \bibinfo{publisher}{Springer Berlin Heidelberg}, pp.
  \bibinfo{pages}{127--138}, \doi{10.1007/978-3-642-02930-1\_11}.

\bibitemdeclare{article}{coecke12}
\bibitem{coecke12}
\bibinfo{author}{B.~\surnamestart Coecke\surnameend} \&
  \bibinfo{author}{R.~\surnamestart Spekkens\surnameend}
  (\bibinfo{year}{2012}): \emph{\bibinfo{title}{Picturing classical and quantum
  Bayesian inference}}.
\newblock {\sl \bibinfo{journal}{Synthese}}
  \bibinfo{volume}{186}(\bibinfo{number}{3}), pp. \bibinfo{pages}{651--696},
  \doi{10.1007/s11229-011-9917-5}.

\bibitemdeclare{article}{Elitzur93}
\bibitem{Elitzur93}
\bibinfo{author}{A.~\surnamestart Elitzur\surnameend} \&
  \bibinfo{author}{L.~\surnamestart Vaidman\surnameend} (\bibinfo{year}{1993}):
  \emph{\bibinfo{title}{Quantum mechanical interaction-free measurements}}.
\newblock {\sl \bibinfo{journal}{Foundations of Physics}}
  \bibinfo{volume}{23}(\bibinfo{number}{7}), pp. \bibinfo{pages}{987--997},
  \doi{10.1007/BF00736012}.

\bibitemdeclare{incollection}{FurberJ13a}
\bibitem{FurberJ13a}
\bibinfo{author}{R.~\surnamestart Furber\surnameend} \&
  \bibinfo{author}{B.~\surnamestart Jacobs\surnameend} (\bibinfo{year}{2013}):
  \emph{\bibinfo{title}{{From Kleisli Categories to Commutative C$^*$-Algebras:
  Probabilistic Gelfand Duality}}}.
\newblock In \bibinfo{editor}{R.~\surnamestart Heckel\surnameend} \&
  \bibinfo{editor}{S.~\surnamestart Milius\surnameend}, editors: {\sl
  \bibinfo{booktitle}{Algebra and Coalgebra in Computer Science}}, {\sl
  \bibinfo{series}{Lecture Notes in Computer Science}} \bibinfo{volume}{8089},
  \bibinfo{publisher}{Springer Berlin Heidelberg}, pp.
  \bibinfo{pages}{141--157}, \doi{10.1007/978-3-642-40206-7\_12}.

\bibitemdeclare{article}{Jacobs13b}
\bibitem{Jacobs13b}
\bibinfo{author}{B.~\surnamestart Jacobs\surnameend} (\bibinfo{year}{2013}):
  \emph{\bibinfo{title}{{On Block Structures in Quantum Computation}}}.
\newblock {\sl \bibinfo{journal}{Electronic Notes in Theoretical Computer
  Science}} \bibinfo{volume}{298}(\bibinfo{number}{0}), pp. \bibinfo{pages}{233
  -- 255}, \doi{10.1016/j.entcs.2013.09.016}.
\newblock \bibinfo{note}{Proceedings of the Twenty-ninth Conference on the
  Mathematical Foundations of Programming Semantics, \{MFPS\} \{XXIX\}}.

\bibitemdeclare{article}{newdirections}
\bibitem{newdirections}
\bibinfo{author}{B.~\surnamestart Jacobs\surnameend} (\bibinfo{year}{2015}):
  \emph{\bibinfo{title}{{New Directions in Categorical Logic for Classical,
  Probabilistic and Quantum Logic}}}.
\newblock {\sl \bibinfo{journal}{Logical Methods in Computer Science}}
  \bibinfo{volume}{11}(\bibinfo{number}{3}).

\bibitemdeclare{inproceedings}{JacobsM12b}
\bibitem{JacobsM12b}
\bibinfo{author}{B.~\surnamestart Jacobs\surnameend} \&
  \bibinfo{author}{J.~\surnamestart Mandemaker\surnameend}
  (\bibinfo{year}{2012}): \emph{\bibinfo{title}{The Expectation Monad in
  Quantum Foundations}}.
\newblock In \bibinfo{editor}{B.~\surnamestart Jacobs\surnameend},
  \bibinfo{editor}{P.~\surnamestart Selinger\surnameend} \&
  \bibinfo{editor}{B.~\surnamestart Spitters\surnameend}, editors: {\sl
  \bibinfo{booktitle}{Quantum Physics and Logic (QPL) 2011}}, {\sl
  \bibinfo{series}{Elect. Proc. in Theor. Comp. Sci.}}~\bibinfo{volume}{95},
  pp. \bibinfo{pages}{143--182}, \doi{10.4204/EPTCS.95.12}.

\bibitemdeclare{incollection}{JacobsM12c}
\bibitem{JacobsM12c}
\bibinfo{author}{B.~\surnamestart Jacobs\surnameend} \&
  \bibinfo{author}{J.~\surnamestart Mandemaker\surnameend}
  (\bibinfo{year}{2016, to appear}): \emph{\bibinfo{title}{Relating Operator
  Spaces via Adjunctions}}.
\newblock In \bibinfo{editor}{J.~Chubb \surnamestart Reimann\surnameend},
  \bibinfo{editor}{V.~\surnamestart Harizanov\surnameend} \&
  \bibinfo{editor}{A.~\surnamestart Eskandarian\surnameend}, editors: {\sl
  \bibinfo{booktitle}{Logic and Algebraic Structures in Quantum Computing and
  Information}}, \bibinfo{series}{Lect. Notes in Logic},
  \bibinfo{publisher}{Cambridge Univ. Press}.
\newblock \urlprefix\url{http://arxiv.org/abs/1201.1272}.

\bibitemdeclare{article}{leifer12}
\bibitem{leifer12}
\bibinfo{author}{M.~\surnamestart Leifer\surnameend} \&
  \bibinfo{author}{R.~\surnamestart Spekkens\surnameend}
  (\bibinfo{year}{2013}): \emph{\bibinfo{title}{Towards a formulation of
  quantum theory as a causally neutral theory of {Bayesian} inference}}.
\newblock {\sl \bibinfo{journal}{Phys. Rev. A}} \bibinfo{volume}{88(5)}, p.
  \bibinfo{pages}{052130}.

\bibitemdeclare{article}{RussoDye66}
\bibitem{RussoDye66}
\bibinfo{author}{B.~\surnamestart Russo\surnameend} \& \bibinfo{author}{H.A.
  \surnamestart Dye\surnameend} (\bibinfo{year}{1966}):
  \emph{\bibinfo{title}{{A Note on Unitary Operators in $C\sp *$-algebras.}}}
\newblock {\sl \bibinfo{journal}{Duke Math. J.}} \bibinfo{volume}{33}, pp.
  \bibinfo{pages}{413--416}, \doi{10.1215/S0012-7094-66-03346-1}.

\bibitemdeclare{book}{TakesakiI}
\bibitem{TakesakiI}
\bibinfo{author}{Masamichi \surnamestart Takesaki\surnameend}
  (\bibinfo{year}{2002}): \emph{\bibinfo{title}{{Theory of Operator Algebra,
  volume I}}}.
\newblock {\sl \bibinfo{series}{Encyclopedia of Mathematical Sciences}}
  \bibinfo{volume}{124}, \bibinfo{publisher}{Springer Verlag},
  \doi{10.1007/978-1-4612-6188-9}.

\bibitemdeclare{book}{TakesakiIII}
\bibitem{TakesakiIII}
\bibinfo{author}{Masamichi \surnamestart Takesaki\surnameend}
  (\bibinfo{year}{2003}): \emph{\bibinfo{title}{{Theory of Operator Algebras,
  volume III}}}.
\newblock {\sl \bibinfo{series}{Encyclopedia of Mathematical Sciences}}
  \bibinfo{volume}{127}, \bibinfo{publisher}{Springer Verlag},
  \doi{10.1007/978-3-662-10453-8}.

\bibitemdeclare{article}{tomiyama}
\bibitem{tomiyama}
\bibinfo{author}{Jun \surnamestart Tomiyama\surnameend} (\bibinfo{year}{1957}):
  \emph{\bibinfo{title}{{On the Projection of Norm One in W*-algebras}}}.
\newblock {\sl \bibinfo{journal}{Proceedings of the Japan Academy}}
  \bibinfo{volume}{33}(\bibinfo{number}{10}), pp. \bibinfo{pages}{608--612},
  \doi{10.3792/pja/1195524885}.

\end{thebibliography}

\end{document}